\documentclass[12pt]{article}
\usepackage{amsmath}
\usepackage{amssymb}
\usepackage{amsthm}
\usepackage{enumitem}
\usepackage{mathrsfs}
\usepackage[dvipdfmx]{graphicx}
\usepackage{float}
\usepackage{bbm}
\usepackage{bm}
\usepackage{comment}
\usepackage{xcolor}
\usepackage[top=3cm,bottom=3cm,right=2.5cm,left=2.5cm]{geometry}

\usepackage[pagebackref]{hyperref}
\hypersetup{
   colorlinks,
    linkcolor={red!60!black},
    citecolor={blue!60!black},
    urlcolor={blue!90!black}
}

\newcommand{\Z}{\mathbb{Z}}
\newcommand{\R}{\mathbb{R}}
\newcommand{\N}{\mathbb{N}}
\newcommand{\Q}{\mathbb{Q}}
\newcommand{\C}{\mathbb{C}}

\newcommand{\SO}{{\rm SO}}

\newcommand{\D}{\mathcal{D}}
\newcommand{\supp}{{\rm supp}}

\newcommand{\Prob}{\mathcal{P}}

\newcommand{\spn}{{\rm span}}

\newcommand{\e}{\varepsilon}
\newcommand{\bigast}{\mathop{\mathchoice
    {\vcenter{\hbox{\Large$\ast$}}}
    {\vcenter{\hbox{\large$\ast$}}}
    {\ast}
    {\ast}}\displaylimits}

\newcommand {\norm}[1] {\left\| {#1} \right\|}

\DeclareMathOperator{\ubdim}{\overline{\dim}_B}

\DeclareMathOperator{\hdim}{\dim_H}
\DeclareMathOperator{\pdim}{\dim_P}

\numberwithin{equation}{section}

\newtheorem{prop}{Proposition}[section]
\newtheorem{thm}[prop]{Theorem}
\newtheorem{cor}[prop]{Corollary}
\newtheorem{lem}[prop]{Lemma}

\theoremstyle{definition}
\newtheorem{dfn}{Definition}[section]

\newtheorem{rem}[prop]{Remark}

\title{Markov-driven self-similar measures and non-homogeneous slicing}

\author{Pablo Shmerkin and Shunsuke Usuki}

\begin{document}
	
\maketitle

\begin{abstract}
We introduce the concept of Markov-driven self-similar measures, a notion that simultaneously generalizes classical self-similar measures, self-conformal measures, dynamically driven self-similar measures and many random measures. This class is also closed under Cartesian products and linear projections. We initiate the study of this class of measures and investigate the $L^q$ dimension of Markov-driven self-similar measures under certain Diophantine conditions. As an application, we obtain the first Furstenberg slicing type results for intersections of Cantor sets in non-homogeneous settings.
\end{abstract}

\section{Introduction and main results}

\subsection{Slicing in non-homogeneous settings}

Furstenberg's slicing conjecture \cite{Fur70} predicted that if $A_1, A_2\subset \R/\Z$ are closed sets invariant under $\times p_1$ and $\times p_2$ respectively, with $p_1,p_2\in\N_{\ge 2}$ multiplicatively independent, then 
\begin{equation}\label{eq:Furstenberg-slicing-conjecture}
    \hdim\bigl((A_1\times A_2)\cap \ell\bigr)\le \max\{\hdim(A_1\times A_2)-1,0\}
\end{equation}
for all lines $\ell$ not parallel to the axes. This was proved, independently, by  M.~Wu \cite{Wu19} and the first author \cite{Shm19}. A simpler proof was later given by T.~Austin \cite{Aus22}. The methods of all these papers work more generally under the assumption that $A_i$ are self-similar sets for \emph{homogeneous} iterated function systems $\{ \lambda_i x + t_j^{i}\}_{j=1}^{k_i}$, satisfying the strong separation condition, with $\lambda_1,\lambda_2\in (0,1)$ multiplicatively independent. This framework was extended to a higher dimensional setting by E.~Corso and the first author \cite{CS24}.

In all these works, and despite the diversity of the approaches, the homogeneity of the underlying iterated function system plays a critical role. The proofs make essential use of the underlying dynamics given by rotation by $\log \lambda_1/\log\lambda_2$ on the circle; the fact that Lebesgue measure is the unique invariant measure for this rotation is used in a crucial way.

In this paper, we obtain the first Furstenberg-type slicing results in non-homogeneous settings. The following is our main theorem in this direction. See  \S\ref{S:Markov-model-self-similar} for the definition of self-similar measures and  \S\ref{S:worst-case-spectrum} for the definition of the $L^q$ dimension $D_q$ of a measure.
\begin{thm} \label{thm:intro-slicing}
    Consider two iterated function systems $\Phi=\{ \lambda_i x + a_i\}_{i\in\mathcal{I}}$, $\Psi=\{\lambda'_j x + a'_j\}_{j\in\mathcal{J}}$ on the real line. Assume that both of them satisfy the open set condition. Further, assume that the following Diophantine condition holds:   there are $K\in\N$, $c_0>0$, and infinitely many $N$ such that for all integers $0\le n_i, m_j\le N\ (i\in\mathcal{I},j\in\mathcal{J})$, with at least one of them non-zero, we have
    \[
        \left| \sum_{i\in\mathcal{I}} n_i\log \lambda_i -\sum_{j\in\mathcal{J}} m_j\log \lambda'_j \right| \ge c_0 N^{-K}.
    \]
    Let $\mu$ and $\nu$ be self-similar measures associated to $\Phi$ and $\Psi$, respectively, and the probability vectors $(p_i)_{i\in\mathcal{I}}$ and $(p'_j)_{j\in\mathcal{J}}$. Let $A=\supp(\mu)$, $B=\supp(\nu)$ be the corresponding self-similar sets.
    
    Then, for all $q\in (1,\infty)$,
    \begin{equation} \label{eq:intro-Lq-convolution}
        D_q(\mu*S_r\nu) = \min\left\{\frac{\tau_{\Phi,p}(q)+\tau_{\Psi,p'}(q)}{q-1},1\right\}\quad\text{for all }r\in\R\setminus\{0\},
    \end{equation}
    where $S_r(x)=rx$ scales by $r$, $\tau_{\Phi,p}(q)$ is the unique real number satisfying $\sum_{i\in\mathcal{I}}p_i^q\lambda_i^{-\tau_{\Phi,p}(q)}=1$ and $\tau_{\Psi,p'}(q)$ is defined similarly. Consequently, for all lines $\ell$ not parallel to the axes,
    \begin{equation} \label{eq:intro-slicing}
        \ubdim\bigl((A\times B)\cap \ell\bigr)\le \max\{\hdim(A\times B)-1,0\}.
    \end{equation}
\end{thm}
We make some remarks on this statement.
\begin{rem} \label{rem:slices}
  Equation \eqref{eq:intro-slicing} follows by a simple observation relating slices to $L^q$ dimensions as $q\to\infty$, see \cite[Lemma 1.8]{Shm19}. More precisely, this is applied to the projections of $\mu\times\nu$, where $\mu,\nu$ are the natural measures on $A,B$, respectively with weights $\lambda_i^{\hdim(A)}$, $(\lambda'_j)^{\hdim(B)}$. 
\end{rem}
\begin{rem}
    Feng, Huang and Rao \cite[Conjecture 1.2]{FHR14} conjectured that if $A,B$ are totally disconnected self-similar sets with contraction ratios $\lambda_i$, $\lambda'_j$, respectively, and $A$ embeds affinely into $B$, then each $\log\lambda_i$ belongs to $\spn_{\Q}\{\log\lambda'_j:j\in\mathcal{J}\}$.

    Slicing results yield non-embeddability results. Under the assumptions of the theorem, letting $A=\supp(\mu)$ and $B=\supp(\nu)$, we have that if $\hdim(A)>0$ and $\hdim(B)<1$, then $A$ cannot be affinely embedded into $B$, for otherwise the graph of the embedding would be a line $\ell$ with $\hdim((A\times B)\cap\ell)=\hdim(A)$, contradicting the slicing bound.
    
    The most general non-embeddability results for linear self-similar sets are due to A.~Algom, M.~Hochman and M.~Wu \cite{AHW24}. Our Diophantine condition is the same as \cite[Condition (d) in Definition 6.2]{AHW24} with $\gamma_i$ consisting of $\log\lambda_i$ and $-\log\lambda'_j$. As a result, Theorem \ref{thm:intro-slicing} recovers and significantly strengthens the results of \cite{AHW24}. Interestingly, even though the Diophantine condition is the same, the methods of this paper are very different to those in \cite{AHW24}.
\end{rem}

\begin{rem}
    See \S\ref{sec:diophantine-convolution} for a discussion of the Diophantine condition. For homogeneous systems with contractions $\lambda,\lambda'$, the condition follows from an elementary argument whenever $\log \lambda/\log\lambda'\notin\mathbb{Q}$. Thus Theorem \ref{thm:intro-slicing} generalizes the results of \cite{Wu19,Shm19,Aus22}.
    It is plausible that \eqref{eq:intro-slicing} holds assuming only that
    \begin{equation} \label{eq:no-exact-overlaps-slicing}
        \prod_{i\in\mathcal I}\lambda_i^{n_i}
        \ne\prod_{j\in\mathcal J}(\lambda'_j)^{m_j}
    \end{equation}
    for all $\mathbf n\in\N_0^{\mathcal I}$ and $\mathbf m\in\N_0^{\mathcal J}$ with $(\mathbf n,\mathbf m)\ne0$. Equivalently, no nonempty compositions from the two systems have the same contraction ratio.
\end{rem}

In a similar spirit, we have the following $L^q$ dimension and slicing result for planar self-similar sets and measures. Recall that $\Pi_{V}$ denotes the orthogonal projection onto the line $V$ in $\R^2$.
\begin{thm} \label{thm:intro-ssm-slices}
    Let $\Phi=\{ \lambda_i R_{\theta_i} x + a_i\}_{i\in\mathcal{I}}$ be an iterated function system on $\R^2$ satisfying the open set condition, where $\lambda_i\in (0,1)$, $R_{\theta_i}$ is rotation by $2\pi\theta_i$, with $\theta_i\in\R/\Z$ and $a_i\in\R^2$. Let $\mu$ be the self-similar measure associated to $\Phi$ and a probability vector $(p_i)_{i\in\mathcal{I}}$, and let $A=\supp(\mu)$. 
    
    Assume that the following Diophantine condition holds: there are $K\in\N$, $c_0>0$ and infinitely many $N$ such that for all integers $0\le n_i\le N\ (i\in\mathcal{I})$, with at least one of them non-zero, we have
    \[
     \left\|\sum_{i\in\mathcal I}n_i\theta_i\right\|_{\R/\Z}
    \ge c_0N^{-K}.
    \]
    Then,
    \begin{equation*}
    D_q(\mu)=\min\left\{\frac{\tau_{\Phi,p}(q)}{q-1},2\right\}
    \end{equation*}
    and
    \[
        D_q(\Pi_{V}\mu) = \min\left\{\frac{\tau_{\Phi,p}(q)}{q-1},1\right\}
    \]
    for all $q\in (1,\infty)$ and one-dimensional linear subspaces $V\subset\R^2$,
    where $\tau_{\Phi,p}(q)$ is the unique real number satisfying $\sum_{i\in\mathcal{I}}p_i^q\lambda_i^{-\tau_{\Phi,p}(q)}=1$.
    In particular,
    \[
        \ubdim(A\cap \ell) \le \max\{\hdim(A)-1,0\}\quad\text{for all lines }\ell\subset \R^2.
    \]
\end{thm}
This result extends \cite[Corollary 8.3]{Shm19} and \cite[Theorem 1.6]{Wu19}, which only considered homogeneous systems. 

\subsection{Markov-driven self-similar measures}
\label{S:intro-Markov-SSM}

The main conceptual innovation of the paper is the introduction of the class of \emph{Markov-driven self-similar measures}. This is a generalization of the \emph{dynamically driven self-similar measures} introduced by the first author in \cite{Shm19}. Let $X$ be a compact metric space and $T:X\to X$ a Borel map, let $\lambda\in (0,1)$ be a contraction parameter, and let $\Delta:X\to \mathcal{A}_d$ be a measurable map to the set of purely atomic Borel probability measures $\mathcal{A}_d$ on the unit ball in $\R^d$. The corresponding family of dynamically driven self-similar measures is defined by
\begin{equation} \label{eq:dyn-ssm-conv}
\mu_x = \bigast_{n=0}^{\infty} S_{\lambda^n}\Delta(T^nx).
\end{equation}
This class of measures includes classical \emph{homogeneous} self-similar measures (corresponding to the case in which $X$ is a singleton). It has the nice property of being invariant under projections, a property that is crucial in the higher dimensional extensions in \cite{CS24}. 

The convolution structure \eqref{eq:dyn-ssm-conv} of dynamically driven self-similar measures is crucial in \cite{Shm19,CS24}. Another key underlying assumption in those works is that $(X,T)$ is uniquely ergodic. In this article, we introduce a more general setting in which the underlying dynamics is given by a Markov kernel, rather than a uniquely ergodic dynamical system. This leads to a substantially more general class of measures with a Markov self-similarity property.

Let $P$ be a Borel Markov kernel on a nonempty Borel subset $X$ of a compact metric space, equipped with the inherited metric (see Section \ref{S:Markov-kernels-operators} for the rigorous definition). For $n\in\N_0$, we denote the law of the orbit started at $x\in X$ up to time $n$ by $P_x^n\in\Prob(X^{n+1})$, and the law of the Markov chain started at $x$ and observed at time $n$ by $Q_x^n\in\Prob(X)$, with $P_x^0=Q_x^0=\delta_x$. We assume that we are given the following data:
\begin{enumerate}
\renewcommand{\labelenumi}{\rm{(\roman{enumi})}}
\item a Borel measurable map $\pi: X\to \mathrm B^d(0,R_0)\ (R_0<\infty)$;
\item a Borel measurable map $\lambda: X\to(0,1)$ such that $0<\lambda_{\min}\le \lambda_{\max}<1$, where
\[
\lambda_{\min}=\inf_{x\in X}\lambda(x),\quad \lambda_{\max}=\sup_{x\in X}\lambda(x).
\]
\end{enumerate}
We introduce the following notation:
\[
    \lambda(x_0^n) = \prod_{i=0}^{n-1} \lambda(x_i) \quad (x_0^n=(x_0,\ldots,x_n)\in X^{n+1}), 
\]
with the convention that an empty product is $1$. Then, for each $n\in\N_0$, we define the Borel measurable map $\pi_n: X^{n+1}\to\R^d$ by
\begin{equation*}
\pi_n(x_0^n)=\sum_{k=0}^n \lambda(x_0^k) \cdot\pi(x_k)\quad (x_0^n=(x_0,\dots,x_n)\in X^{n+1}).
\end{equation*}
We use the same notation on $X^{\mathbb{N}_0}$, writing
\begin{equation*}
\pi_n(x_0^\infty)=\sum_{k=0}^n \lambda(x_0^k)\cdot \pi(x_k)\quad (x_0^\infty=(x_0,x_1,\dots)\in X^{\mathbb{N}_0}).
\end{equation*}
For $x\in X$ and $n\in\N_0$, we define a Borel probability measure $\mu_{x,n}$ on $\R^d$ by
\begin{equation*}
\mu_{x,n}=\pi_nP_x^n.
\end{equation*}
Write $R=(1-\lambda_{\max})^{-1}$, and note that $\supp(\mu_{x,n})\subset \mathrm B^d(0,RR_0)$ for all $x\in X$ and $n\in\N_0$.  Furthermore, for $x\in X$ and $n,m\in\N$,
\begin{equation} \label{eq:time-n+m-decomposition}
 \mu_{x,n+m} = \int_{X^{n+1}} \tau_{\pi_{n-1}(x_0^{n-1})} S_{\lambda(x_0^{n})} \mu_{x_{n},m}  \ d P_x^n(x_0^n),
\end{equation}
where $\tau_a(y) = y + a$. Indeed, this follows since
\[
\pi_{n+m}(x_0^{n+m}) = \pi_{n-1}(x_0^{n-1}) + \lambda(x_0^n) \cdot\pi_m(x_n^{n+m}).
\]

For every one-sided orbit $x_0^\infty\in X^{\mathbb{N}_0}$, the series
\[
\pi_\infty(x_0^\infty)=\sum_{k=0}^\infty \lambda(x_0^k)\pi(x_k)
\]
converges absolutely, with $|\pi_\infty-\pi_n|\le RR_0\lambda_{\max}^{n+1}$ uniformly on $X^{\mathbb N_0}$. Hence $\mu_{x,n}$ converges weak*-ly to $\pi_\infty P_x^\infty$, where $P_x^\infty\in\mathcal{P}(X^{\N_0})$ is the law of the full one-sided Markov orbits started from $x$.

\begin{dfn}[Markov-driven self-similar measure]\label{Markov_self_similar_measures}
Given $X, P, \pi, \lambda$ as above, we refer to the tuple
\[
 \mathcal{X} = (X, P, \pi, \lambda)
\]
as a \emph{Markov self-similar model}. For $x\in X$, we write
\begin{equation*}
\mu_x=\lim_{n\to\infty}\mu_{x,n}
\end{equation*}
and call the family $(\mu_x)_{x\in X}$ the family of {\it Markov-driven self-similar measures} generated by $\mathcal{X}$, or simply {\it Markov self-similar measure} if the context is clear.
\end{dfn}
Thus
\[
\mu_{x}=\pi_\infty P_{x}^\infty=\int_{X^{\mathbb{N}_0}}\delta_{\pi_\infty(x_0^\infty)}\,dP_{x}^\infty(x_0^\infty).
\]

Taking $m\to\infty$ in \eqref{eq:time-n+m-decomposition}, we get
\begin{equation}\label{self_similarity_of_Markov_self_similar_measures}
\mu_{x}=\int_{X^{n+1}} \tau_{\pi_{n-1}(x_0^{n-1})} S_{\lambda(x_0^n)} \mu_{x_n} \, dP_{x}^n(x_0^n).
\end{equation}
This is the Markov self-similarity relation: the measure $\mu_x$ is obtained by averaging scaled-down and translated copies of the measures $\mu_{x_n}$, where $x_n$ is distributed according to $Q_{x_0}^n$.

We note the simple fact that the Markov self-similarity relation \eqref{self_similarity_of_Markov_self_similar_measures} characterizes the family $(\mu_x)_{x\in X}$. We will use this fact repeatedly, without further reference, when we provide examples of Markov self-similar measures in Section \ref{sec:examples}. 
\begin{lem}\label{lem:Markov-self-similarity-uniqueness}
Given a Markov self-similar model $\mathcal{X}$.
Then, the family $(\mu_x)_{x\in X}$ of the Markov-self-similar measures is the unique Borel family of probability measures with uniformly bounded supports satisfying \eqref{self_similarity_of_Markov_self_similar_measures} for $n=1$.
\end{lem}
\begin{proof}
Iterating \eqref{self_similarity_of_Markov_self_similar_measures} with $n=1$ and using $\lambda(x_0^n)\le\lambda_{\max}^n\to0$ shows that the measure at $x$ equals $\lim_{n\to\infty}\mu_{x,n}=\mu_x$.
\end{proof}


We discuss examples of Markov self-similar measures in Section \ref{sec:examples}. For now, we note that this class is invariant under Cartesian products and linear projections (hence, also under convolutions), and includes self-similar measures, self-conformal measures, dynamically driven self-similar measures, and many random constructions.

Our main result on Markov self-similar measures is Theorem \ref{thm:diophantine-transversal-Lq}, which provides a lower bound for the $L^q$ dimension of Markov self-similar measures under certain assumptions. Applying the theorem to the convolution and planar models yields Theorems \ref{thm:intro-slicing} and \ref{thm:intro-ssm-slices}.

The assumptions of Theorem \ref{thm:diophantine-transversal-Lq} are often hard to verify; even when they hold, the provided bound may not be sharp. Still, we think that it illustrates the power of the Markov self-similar framework, and we hope other applications will be found in the future.

\subsection{Overview of the proofs}

Our approach borrows many ideas from \cite{Shm19,CS24}, but there are also substantial differences. The approach of \cite{Shm19,CS24} is based on a sub-multiplicative relation among the $L^q$ norms of the measures $\mu_{x,n}$ generated by a dynamically driven self-similar model. This sub-multiplicativity, together with versions of the subadditive ergodic theorem for uniquely ergodic systems, implies that the ``worst case'' $L^q$ dimension of the measures $\mu_{x,n}$ is also the ``typical case''.

To our knowledge, there is no suitable version of the sub-multiplicative relation for Markov self-similar measures. Instead, we introduce a new approach based on the \emph{worst case spectrum}, introduced in \S\ref{S:worst-case-spectrum}. This yields trivially a lower bound for the $L^q$ spectrum of all the measures $\mu_x$ generated by the model. Perhaps surprisingly, a lot of the machinery of \cite{Shm19,CS24} can be adapted to the worst case spectrum, except that we work with finitary statements instead of taking limits. The Diophantine assumptions, together with a transversality condition introduced in Definition \ref{dfn:transversal}, are used as a finitary replacement of the crucial exponential separation property from \cite{Shm19, CS24}.

\subsection{AI statement}

All the main ideas in the paper, including the concept of Markov self-similar measures, and the approach to $L^q$ dimensions via the worst case spectrum, were obtained by the authors in 2025 and early 2026 without any AI assistance. The authors used LLMs for literature review, to draft proofs of some technical statements, and for checking mathematical correctness. The authors have verified and substantially rewritten all AI-generated content. 

\subsection{Acknowledgements}

This material is partially based upon work supported by the Swedish Research Council under grant no. 2021-06594 while the authors were in residence at Institut Mittag-Leffler in Djursholm, Sweden, during the semester \emph{Interactions between fractal geometry, harmonic analysis, and dynamics in 2026}.

P.S. was partially supported by a Discovery grant from the Natural Sciences and Engineering Research Council of Canada (NSERC). S.U. was partially supported by a KAKENHI grant from the Japan Society for the Promotion of Science (JSPS) (Grant Number 26K17002).

\section{Markov-driven self-similar measures}\label{sec:Markov_self_similar_measures}

\subsection{Notation}

We denote $\N=\{1,2,\ldots\}$ and $\N_0=\{0,1,2,\ldots\}$. All logarithms are to base $2$. We write $f\nu$ for the push-forward of a measure $\nu$ under a Borel map $f$. 

For nonnegative quantities $a,b$, we write $a\lesssim b$ or $b\gtrsim a$ if $a\le Cb$ for a constant $C>0$, and $a\sim b$ if $a\lesssim b$ and $b\lesssim a$. Implicit constants may depend on the fixed model and ambient dimension; additional dependencies are indicated by subscripts.

For $u=(u_1,\dots,u_d)\in\R^d$, we write $|u|$ for the Euclidean norm of $u$ and $|u|_\infty=\max_{1\leq i\leq d}|u_i|$. For $r>0$, $B(u,r)=\{v\in\R^d:|v-u|\leq r\}$ denotes the closed ball with center $u$ and radius $r$.

For $m\in\N_0$, let $\D_m=\{2^{-m}(k+[0,1)^d):k\in\Z^d\}$ be the dyadic partition of $\R^d$. For $E\subset\R^d$, write $\D_m(E)=\{Q\in\D_m:Q\cap E\ne\varnothing\}$; for $x\in\R^d$, $\D_m(x)$ denotes the unique cube containing $x$.

The Grassmannian of linear hyperplanes in $\R^d$ is denoted by $\mathbb{G}(d,d-1)$. For $V\in\mathbb{G}(d,d-1)$, we write $\Pi_V:\R^d\to V$ for the orthogonal projection onto $V$. Dyadic partitions and discretizations on $V\in \mathbb{G}(d,d-1)$ are defined using fixed orthonormal coordinates.

For a cube $Q$ and $a>0$, write $aQ$ for the concentric cube with side length $a$ times that of $Q$.

\subsection{Borel Markov kernels and Markov operators}
\label{S:Markov-kernels-operators}

We recall some basic definitions and properties of Borel Markov kernels and Markov operators, and introduce notation along the way.

Throughout the paper, a phase space $X$ is a nonempty Borel subset of a compact metric space, with the inherited metric and Borel $\sigma$-algebra. In particular, $X$ is a standard Borel space. Kernels are defined on $X$. Supports of measures on $X$ are taken in its relative topology, and covering numbers use the inherited metric.
We write $\mathcal{B}_b(X)$ for the space of all bounded Borel measurable functions $f:X\to\C$.
We also write $\Prob(X)$ for the set of all Borel probability measures on $X$.

A {\it Borel Markov kernel} on $X$ is a map $P:X\to\Prob(X)$ such that, for every Borel set $E\subset X$, the function $x\mapsto P(x)(E)$ is Borel measurable. Equivalently, for every $f\in\mathcal{B}_b(X)$, the function
$x\mapsto \int_X f\,dP(x)$ is Borel measurable. For $x\in X$, we often write $P_x=P(x)\in\Prob(X)$.

Given a sequence $x=(x_0,\ldots,x_n)$, we write $x_i^j=(x_i,x_{i+1},\ldots,x_j)$ for $0\leq i\leq j\leq n$. Sometimes we also write $x=x_0^n$ to emphasize the length of the sequence $x$.

For $x\in X$ and $n\in\N$, we write $P_x^n\in\Prob(X^{n+1})$ for the distribution of orbits starting from $x$ up to time $n$, that is, $P_x^n$ is the Borel probability measure on $X^{n+1}$ defined by
\begin{equation*}
\int_{X^{n+1}} f\,dP_x^n=\int_X\int_X\cdots\int_X\int_X f(x_0^n)\,dP_{x_{n-1}}(x_n)dP_{x_{n-2}}(x_{n-1})\cdots dP_{x_1}(x_2)dP_x(x_1)
\end{equation*}
for every $f\in\mathcal{B}_b(X^{n+1})$. We also denote by $Q_x^n$ the distribution of the Markov chain started at $x$ and observed at time $n$, that is, $Q_x^n$ is the Borel probability measure on $X$ defined by
\begin{equation*}
    \int_X f\,dQ_x^n= \int_{X^{n+1}} f(x_n)\,dP_x^n(x,x_1,\ldots,x_n)
\end{equation*}
for every bounded Borel measurable function $f\in\mathcal{B}_b(X)$. We set $P_x^0=Q_x^0=\delta_x$.

By Kolmogorov's extension theorem, for each $x\in X$ there is a unique Borel probability measure $P_x^\infty$ on $X^{\mathbb{N}_0}$ whose marginal on coordinates $0,\ldots,n$ is $P_x^n$ for every $n\in\N_0$.
This is the law of the full one-sided Markov orbit started from $x$.

The maps $x\mapsto P_x^n$, $x\mapsto Q_x^n$ and $x\mapsto P_x^\infty$ are Borel kernels. We call a finite or one-sided sequence $(x_i)$ a \emph{Markov orbit} if $x_{i+1}\in\supp(P_{x_i})$ at every step. The relation $\{(x,y):y\in\supp(P_x)\}$ is Borel, and $P_x^\infty$-almost every sequence is a Markov orbit.

\subsection{Support renormalization}

\label{S:support-renormalization}

From now, we will work with a Markov self-similar model $\mathcal{X}=(X,P,\pi,\lambda)$ as in Section \ref{S:intro-Markov-SSM}.

Let $H(y)=ry+b$ for some $r>0$ and $b\in\R^d$, and define
\[
\pi'(x)=r\pi(x)+(1-\lambda(x))b.
\]
Let $(\mu'_x)_{x\in X}$ be the family generated by the model $(X,P,\pi',\lambda)$. Then $\mu'_x=H\mu_x$ for all $x\in X$. By choosing $H$ appropriately, we can arrange that
\[
\supp(\mu'_x)\subset [0,1)^d
\qquad\text{for all }x\in X.
\]
For simplicity, the estimates below are written under this normalization; the same proofs apply to supports in $\mathrm B^d(0,RR_0)$, with constants depending additionally on $RR_0$. 

\subsection{Stopped Markov self-similar measures}

In general, for fixed $x\in X$ and large $n$, the contraction ratios $\lambda(x_0^n)$ can vary exponentially among different orbits started from $x$. In practice, it is useful to work with discrete approximations that live at essentially one dyadic scale. We can achieve this by stopping the Markov chain at the first time the contraction ratio crosses a certain threshold. We introduce some notation to formalize this idea.

\begin{dfn}\label{dfn:stopped-path-law}
Fix $t\ge0$. For a one-sided orbit $x_0^\infty\in X^{\N_0}$, define the stopping time
\[
T_t(x_0^\infty)=\min\{n\ge 1:\lambda(x_0^n)<2^{-t}\}.
\]
For $x\in X$, let $P_x^{[t]}$ denote the law of the stopped orbit segment $x_0^{T_t}$ under the Markov chain started from $x$, viewed as a measure on $\bigsqcup_{n=1}^{\infty}X^{n+1}$. We write
\[
    \widetilde{\pi}:\bigsqcup_{n=1}^{\infty}X^{n+1}\to\R^d,
    \qquad \widetilde{\pi}(x_0^n)=\pi_n(x_0^n)\quad(x_0^n\in X^{n+1}, n\in\N).
\]
for the universal coding map.
We define the \emph{stopped Markov self-similar measure} by
\begin{equation*}
\mu_{x,[t]}=\widetilde{\pi}P_x^{[t]}.
\end{equation*}
Finally, let $Q_x^{[t]}$ denote the law of the stopped state $x_{T_t}$.
\end{dfn}

The stopped Markov self-similar measure should be seen as an approximation of $\mu_x$ at scale $2^{-t}$.
If we consider the stopped coding map
\begin{equation*}
\pi_{T_t}:X^{\N_0}\to\R^d,\qquad \pi_{T_t}(x_0^\infty)=\pi_{T_t(x_0^\infty)}(x_0^\infty),
\end{equation*}
$\mu_{x,[t]}$ is written as
\begin{equation*}
\mu_{x,[t]}=\pi_{T_t}P^\infty_x.
\end{equation*}

Note that in the definitions of $P_x^{[t]}$, only $O(t+1)$ values of $n$ are relevant, since $\lambda_{\max}<1$ implies that $T_t=O(t+1)$ uniformly in the orbit. We also remark that both $x\mapsto P_x^{[t]}$ and $x\mapsto Q_x^{[t]}$ are Borel kernels. 

For later use, set
\[
\alpha_{\min}=\log(1/\lambda_{\max}),
\qquad
\alpha_{\max}=\log(1/\lambda_{\min}).
\]
The family of stopped approximations enjoys a Markov self-similarity relation that is very similar to \eqref{self_similarity_of_Markov_self_similar_measures}, and this will be a key tool in our analysis.
\begin{prop}\label{prop:stopped-Markov-self-similarity}
Let
\[
E_{s,n}=\{x_0^n\in X^{n+1}: \lambda(x_0^{n-1})\ge 2^{-s}>\lambda(x_0^n)\}
\]
Fix $s\ge0$ and $t\ge\alpha_{\max}$. For $n\ge1$ and $x_0^n\in E_{s,n}$, define 
\[
\Delta_s(x_0^n)=\log\left(\frac{2^{-s}}{\lambda(x_0^n)}\right).
\]
Then $0<\Delta_s(x_0^n)\le \alpha_{\max}$ and, for $x\in X$,
\begin{equation}\label{eq:stopped-self-similarity}
\mu_{x,[s+t]}
=
\sum_{n=1}^{\infty}
\int_{E_{s,n}}
\tau_{\pi_{n-1}(x_0^{n-1})}S_{\lambda(x_0^n)}
\mu_{x_n,[t-\Delta_s(x_0^n)]}
\ dP_{x}^n(x_0^n).
\end{equation}
\end{prop}

\begin{proof}
On $E_{s,n}$, the identity $\lambda(x_0^n)=2^{-s-\Delta_s(x_0^n)}$ yields
\[
T_{s+t}(x_0^\infty)=n+T_{t-\Delta_s(x_0^n)}(x_n^\infty).
\]
Combining this with
\[
\pi_{n+k}(x_0^{n+k})=\pi_{n-1}(x_0^{n-1})+\lambda(x_0^n)\pi_k(x_n^{n+k})
\]
and the Markov property yields \eqref{eq:stopped-self-similarity}.
\end{proof}

\section{Examples of Markov self-similar measures}
\label{sec:examples}

We illustrate the generality of the Markov self-similar measure point of view by showing that it includes many classes of fractal measures that have so far been studied separately in the literature, including standard self-similar measures, self-conformal measures, dynamically driven self-similar measures, many random constructions, and Furstenberg CP-processes. We also show that the class of Markov self-similar measures is closed under orthogonal projections, Cartesian products and convolutions.  We emphasize that our later results on $L^q$ dimensions of Markov self-similar measures apply only to a small subset of the examples below; our goal here is to illustrate the generality of the framework, with a view towards possible future applications.

\subsection{Classical self-similar measures}

\label{S:Markov-model-self-similar}

Consider an iterated function system (IFS) of similarities on $\R^d$ of the form
\[
  f_i(x) = \lambda_i O_i x + t_i, \quad i=1,\ldots,k,
\]
where $\lambda_i\in (0,1)$, $O_i$ is an orthogonal transformation and $t_i\in\R^d$. Let $p=(p_1,\ldots,p_k)$ be a probability vector. Then the self-similar measure $\mu$ associated to this IFS and probability vector is the unique Borel probability measure on $\R^d$ satisfying
\[
  \mu = \sum_{i=1}^k p_i f_i \mu.
\]

We can represent $\mu$ as a Markov self-similar measure as follows. Let $G$ be the closed subgroup of $\mathrm O(d)$ generated by $O_1,\ldots,O_k$. Let $X=G\times\{1,\ldots,k\}$, and define a continuous Markov kernel $P$ on $X$ by
\[
P_{(O,i)} = \sum_{j=1}^k p_j \delta_{(OO_i,j)}.
\]
That is, $(O,i)$ transitions to $(OO_i,j)$ with probability $p_j$. Define $\pi: X\to\R^d$ by $\pi(O,i)=O t_i$, and $\lambda: X\to(0,1)$ by $\lambda(O,i)=\lambda_i$. Then, for any $(O,i)\in X$, the Markov self-similar measure $\mu_{(O,i)}$ generated by $(X,P,\pi,\lambda)$ is equal to $O(f_i\mu)$. 

\subsection{Self-conformal measures}
\label{S:Markov-model-self-conformal}

For simplicity, we work in dimension $1$. Let $\{ f_i\}_{i\in\mathcal{I}}$ be a finite family  $C^{1+\alpha}$ contractions of $[0,1]$ with $f_i'>0$ on $[0,1]$, and let $p=(p_i)_{i\in\mathcal{I}}$ be a probability vector. Then the self-conformal measure $\mu$ associated to this IFS and probability vector is the unique Borel probability measure on $[0,1]$ satisfying
\[
    \mu = \sum_{i\in\mathcal{I}} p_i f_i \mu.
\]

We can embed $\mu$ into a family of Markov self-similar measures by using the limit geometries determined by the past. Let $\Omega^{-}=\mathcal{I}^{-\N}$ be the space of left-infinite words, and for $\omega^-=(\ldots,\omega_{-2},\omega_{-1})\in\Omega^{-}$ let
\[
    \mathbf{i}_n(\omega^-)=\omega_{-n}\cdots \omega_{-1}.
\]
Write $f_{\mathbf{i}_n(\omega^-)}=f_{\omega_{-n}}\circ\cdots\circ f_{\omega_{-1}}$, and let $H_{\mathbf{i}_n(\omega^-)}$ be the positively oriented affine map sending $f_{\mathbf{i}_n(\omega^-)}([0,1])$ onto $[0,1]$. By the classical limit-geometry theorem for $C^{1+\alpha}$ conformal systems (see \cite[Theorems~5.9 and 6.1]{BF97}), the maps
\[
    G_{\omega^-,n}=H_{\mathbf{i}_n(\omega^-)}\circ f_{\mathbf{i}_n(\omega^-)}
\]
converge in $C^1$ to a diffeomorphism
\[
    G_{\omega^-}:[0,1]\to [0,1].
\]
Moreover, the map $\omega^-\mapsto G_{\omega^-}$ is continuous.

For each $\omega^-\in \Omega^{-}$ and $i\in\mathcal I$, let $\omega^- i\in\Omega^-$ denote the past obtained by appending $i$ at the end. Then
\[
    L_{\omega^-,i}=G_{\omega^-}\circ f_i\circ G_{\omega^- i}^{-1}
\]
is affine. Indeed, for every $n$,
\[
    G_{\omega^-,n}\circ f_i\circ G_{\omega^- i,n+1}^{-1}
    =
    H_{\mathbf{i}_n(\omega^-)}\circ H_{\mathbf{i}_n(\omega^-)i}^{-1},
\]
and the right-hand side is affine. Passing to the limit shows that $L_{\omega^-,i}$ is affine. Define $\lambda(\omega^-,i)$ and $\pi(\omega^-,i)$ by
\[
    L_{\omega^-,i}(x)=\lambda(\omega^-,i)x+\pi(\omega^-,i).
\]

We now define a Markov self-similar model. Set $X=\Omega^-\times \mathcal I$,
and define a continuous Markov kernel on $X$ by
\[
    P_{(\omega^-,i)}
    =
    \sum_{j\in\mathcal I} p_j\,\delta_{(\omega^- i,j)}.
\]
In other words, the current symbol $i$ is appended to the past, and then the next symbol is sampled with law $p$. Since
\begin{align*}
    G_{\omega^-}(f_i\mu) &=
    \sum_{j\in\mathcal I} p_j\, G_{\omega^-}(f_i f_j\mu) \\
    &=
    \sum_{j\in\mathcal I} p_j\, \tau_{\pi(\omega^-,i)}S_{\lambda(\omega^-,i)}G_{\omega^- i}(f_j\mu),
\end{align*}
the family $(G_{\omega^-}(f_i\mu))_{(\omega^-,i)\in X}$ satisfies the Markov self-similarity relation.  Thus the self-conformal measure is embedded, up to the limit diffeomorphisms $G_{\omega^-}$, into a family of Markov self-similar measures.

\subsection{One-dimensional projections of self-affine measures}
\label{S:Markov-model-self-affine-projections}

Consider an IFS of affine contractions on $\R^d$ of the form
\[
    f_i(x)=L_i x+t_i,\qquad i=1,\ldots,k,
\]
where each $L_i$ is invertible and $\max_i\norm{L_i}<1$. Let $p=(p_1,\ldots,p_k)$ be a probability vector, and let $\mu$ be the self-affine measure satisfying
\[
    \mu=\sum_{i=1}^k p_i f_i\mu.
\]
For $u\in\mathbb S^{d-1}$, write
\[
    \Pi_u(x)=\langle u,x\rangle,
    \qquad
    T_i u=\frac{L_i^T u}{\norm{L_i^T u}}.
\]
We define a Markov self-similar model on $\R$ with phase space $X=\mathbb S^{d-1}\times\{1,\ldots,k\}$ by
\[
    P_{(u,i)}=\sum_{j=1}^k p_j\delta_{(T_i u,j)},
    \qquad
    \pi(u,i)=\langle u,t_i\rangle,
    \qquad
    \lambda(u,i)=\norm{L_i^T u}.
\]
The contraction function satisfies
\[
    \min_{1\le i\le k}\norm{(L_i^T)^{-1}}^{-1}
    \le \lambda(u,i)
    \le \max_{1\le i\le k}\norm{L_i}<1.
\]
Since
\[
\Pi_u\circ f_i=\tau_{\langle u,t_i\rangle}S_{\|L_i^Tu\|}\circ\Pi_{T_i u},
\]
the family $(\Pi_u(f_i\mu))_{(u,i)\in X}$ satisfies \eqref{self_similarity_of_Markov_self_similar_measures} for $n=1$. Lemma \ref{lem:Markov-self-similarity-uniqueness} therefore yields $\mu_{(u,i)}=\Pi_u(f_i\mu)$.
 In particular, $\Pi_u\mu=\sum_{i=1}^k p_i\mu_{(u,i)}$. Thus every one-dimensional projection of $\mu$ is embedded into the family generated by this Markov model through its first-generation pieces.

Let $\overline{T}_i:\mathbb P^{d-1}\to\mathbb P^{d-1}$ be the induced projective map $\overline{T}_i([u])=[L_i^T u]$. Any probability measure $\eta^*$ on $\mathbb P^{d-1}$ satisfying
\[
    \eta^*=\sum_{i=1}^k p_i\overline{T}_i\eta^*
\]
is a Furstenberg measure for the dual matrix cocycle. The first coordinate of the Markov chain $(u_n,i_n)$ projects to the corresponding Furstenberg random walk on $\mathbb P^{d-1}$.

\subsection{Dynamically driven self-similar measures}

\label{S:Markov-model-dynamically-driven}

Recall the definition of dynamically driven self-similar measures $(\mu_x)_{x\in X}$ from \S\ref{S:intro-Markov-SSM}. We show that they can be seen as Markov self-similar measures. The phase space of the Markov chain is $Y=\bigcup_{x\in X} \{x\}\times\supp(\Delta(x))$, with the metric inherited from $X\times \mathrm B^d(0,1)$. The chain is given by
\[
P_{(x,t)} = \sum_{u\in\supp(\Delta(Tx))} \Delta(Tx)(u) \cdot \delta_{(Tx,u)}.
\] 
In other words, $(x,t)$ jumps to $(Tx,u)$ with probability $\Delta(Tx)(u)$ for each $u\in\supp(\Delta(Tx))$, independently of the previous $t$. Note that the first component of the chain is deterministically given by the dynamical system $(X,T)$. The map $\pi:Y\to\R^d$ is given by $\pi(x,t)=t$, and $\lambda:Y\to(0,1)$ is given by $\lambda(x,t)=\lambda$. It is straightforward to see that the generated measures satisfy
\[
\mu_{(x,t)}=\delta_t*S_\lambda\mu_{Tx},
\qquad
\mu_x=\int\mu_{(x,t)}\,d\Delta(x)(t).
\]

\subsection{Projections of Markov self-similar measures}

\label{S:Markov-model-projections}

A key property of dynamically driven self-similar measures is that they are closed under orthogonal projections. More precisely, given a model $(X,T,\Delta,\lambda)$ for a dynamically driven self-similar measure, and an orthogonal projection $\Pi_V:\R^d\to V$ for some $k$-dimensional subspace $V$ of $\R^d$, the family of measures $(\Pi_V\mu_x)_{x\in X}$ is also a family of dynamically driven self-similar measures, generated by the model $(X,T,x\mapsto\Pi_V\Delta(x),\lambda)$. On the other hand, classical self-similar measures are not closed under projections (unless the orthogonal transformations $O_i$ are trivial). 

Markov self-similar measures are also closed under projections. More precisely, given a Markov model $\mathcal{X}=(X,P,\pi,\lambda)$, and an orthogonal projection $\Pi_V:\R^d\to V$, the family of measures $(\Pi_V\mu_x)_{x\in X}$ is also a family of Markov self-similar measures, generated by the model $\Pi_V\mathcal{X}=(X,P,\Pi_V\circ\pi,\lambda)$. 

In particular, projections of general classical self-similar measures are natural examples of Markov self-similar measures, which only in very special cases are either classical or dynamically driven self-similar measures.

\subsection{Products of Markov self-similar models}

\label{S:Markov-model-products}

Let
\[
\mathcal{X}_i=(X_i,P^{(i)},\pi_i,\lambda_i)
\qquad (i=1,2)
\]
be Markov self-similar models on $\R^{d_i}$, and let $(\mu^{(i)}_x)_{x\in X_i}$ be the corresponding families of Markov self-similar measures. There is a Markov self-similar model $\mathcal{Z}$ on $\R^{d_1+d_2}$ generating the family of measures
\[
\mu^{(1)}_{x_1}\times S_{2^u}\mu^{(2)}_{x_2}
\qquad ((x_1,x_2)\in X_1\times X_2),
\]
for $u$ in a certain bounded interval. This is done by the similar manner as \cite[Lemma 7.1]{Shm19}.

Write $\alpha_i(x_i)=\log(1/\lambda_i(x_i))$ for $x_i\in X_i$ and $i=1,2$.
We assume that $\alpha_1(x_1)<\alpha_2(x_2)$ for every $x_i\in X_i$. This can always be arranged by slowing the first chain as follows. Choose $N\ge2$ with $N^{-1}\sup_{X_1}\alpha_1<\inf_{X_2}\alpha_2$, and replace $X_1$ by $X_1\times\{0,\ldots,N-1\}$. Use contraction $\lambda_1(x)^{1/N}$ at every $(x,k)$, and coding map $\pi_1(x)$ at phase $0$ and $0$ at the other phases. For $k<N-1$, move deterministically to $(x,k+1)$; from $(x,N-1)$, move to $(y,0)$ with $y$ distributed according to $P^{(1)}_x$. The generated measure at $(x,0)$ is exactly $\mu^{(1)}_x$, since each block of $N$ steps contributes $\pi_1(x)$ and has total contraction $\lambda_1(x)$.
\begin{prop}
Assume that
\[
\alpha_1(x_1)<\alpha_2(x_2)
\qquad\text{for all }x_1\in X_1,\ x_2\in X_2.
\]
Set
\[
Z=\{(x_1,x_2,u): x_i\in X_i,\ 0\le u<\alpha_2(x_2)\}.
\]
Define a Markov kernel $R$ on $Z$ by
\[
R_{(x_1,x_2,u)}
=
\begin{cases}
P_{x_1}^{(1)}\times\delta_{x_2}\times\delta_{u+\alpha_1(x_1)},&\text{if }u+\alpha_1(x_1)<\alpha_2(x_2),\\
P_{x_1}^{(1)}\times P_{x_2}^{(2)}\times\delta_{u+\alpha_1(x_1)-\alpha_2(x_2)},&\text{if }u+\alpha_1(x_1)\ge\alpha_2(x_2),
\end{cases}
\]
and define
\[
\rho(x_1,x_2,u)
=
\begin{cases}
(\pi_1(x_1),0), & \text{if }u+\alpha_1(x_1)<\alpha_2(x_2),\\
(\pi_1(x_1),2^u\pi_2(x_2)), & \text{if }u+\alpha_1(x_1)\ge \alpha_2(x_2),
\end{cases}
\]
\[
\Lambda(x_1,x_2,u)=\lambda_1(x_1).
\]
Then
\[
\mathcal Z=(Z,R,\rho,\Lambda)
\]
is a Markov self-similar model, and its generated measures satisfy
\[
\mu^\mathcal{Z}_{(x_1,x_2,u)}=\mu^{(1)}_{x_1}\times S_{2^u}\mu^{(2)}_{x_2}.
\]
\end{prop}
\begin{proof}
Define
\[
\nu_{(x_1,x_2,u)}=\mu^{(1)}_{x_1}\times S_{2^u}\mu^{(2)}_{x_2}
\qquad ((x_1,x_2,u)\in Z).
\]
If
\[
u+\alpha_1(x_1)<\alpha_2(x_2),
\]
then
\[
\nu_{(x_1,x_2,u)}
=
\int_{X_1}\tau_{(\pi_1(x_1),0)}S_{\lambda_1(x_1)}
\nu_{(y_1,x_2,u+\alpha_1(x_1))}\,dP^{(1)}_{x_1}(y_1),
\]
because $\lambda_1(x_1)2^{u+\alpha_1(x_1)}=2^u$.

If
\[
u+\alpha_1(x_1)\ge \alpha_2(x_2),
\]
then
\[
0\le u+\alpha_1(x_1)-\alpha_2(x_2)<\alpha_1(x_1)<\alpha_2(y_2)
\]
for every $y_2\in X_2$, so the next state indeed belongs to $Z$. Using the self-similarity relation for $\mu^{(2)}_{x_2}$, we get
\[
\nu_{(x_1,x_2,u)}
=
\int_{X_1}\int_{X_2}
\tau_{(\pi_1(x_1),2^u\pi_2(x_2))}S_{\lambda_1(x_1)}
\nu_{(y_1,y_2,u+\alpha_1(x_1)-\alpha_2(x_2))}
\,dP^{(2)}_{x_2}(y_2)\,dP^{(1)}_{x_1}(y_1),
\]
because
\[
\lambda_1(x_1)2^{u+\alpha_1(x_1)-\alpha_2(x_2)}
=
2^u\lambda_2(x_2).
\]
Therefore, we showed that $\nu_{(x_1,x_2,u)}$ satisfies (\ref{self_similarity_of_Markov_self_similar_measures}) for $n=1$ in the both cases, and complete the proof.
\end{proof}

Of course, the above construction can be iterated to show that the Cartesian products of finitely many Markov self-similar measures are again Markov self-similar measures.

\subsection{Convolutions of Markov self-similar measures}\label{sec:convolution-model}

Combining the constructions in \S\ref{S:Markov-model-projections} and \S\ref{S:Markov-model-products}, we see that if $\mathcal{X}_i$, $i=1,2$ are Markov self-similar models on the same space $\R^d$, then there is a model $\mathcal{Z}$ generating the family of measures $(\mu^{(1)}_{x_1}*S_{2^u}\mu^{(2)}_{x_2})_{(x_1,x_2)\in X_1\times X_2}$, for $u$ in a certain bounded interval. Explicitly, we only need to change the definition of $\rho$ in \S\ref{S:Markov-model-products} to 
\[
\rho(x_1,x_2,u)
=
\begin{cases}
\pi_1(x_1), & \text{if }u+\alpha_1(x_1)<\alpha_2(x_2),\\
\pi_1(x_1)+2^u\pi_2(x_2), & \text{if }u+\alpha_1(x_1)\ge \alpha_2(x_2),
\end{cases}
\]
In particular, convolutions of self-similar measures on $\R^d$ can be represented as Markov self-similar measures while, again, only in very special cases are they either classical or dynamically driven self-similar measures.

\subsection{Random fractal measures}

Fix an integer $M\ge2$, let $\mathcal I=\{0,\ldots,M-1\}^d$, and let $\mathcal A$ be the family of non-empty subsets of $\mathcal I$. Write $\mathcal I^*=\bigcup_{n\ge0}\mathcal I^n$, including the empty word $\varnothing$, and equip $\mathcal T=\mathcal A^{\mathcal I^*}$ with the product topology. For $T=(D_{\mathbf{i}})_{\mathbf{i}\in\mathcal I^*}\in\mathcal T$, let $\nu_T$ be the probability measure obtained by starting with unit mass on $[0,1]^d$ and, at every selected cube indexed by $\mathbf{i}$, splitting its mass equally among the $M$-adic subcubes indexed by $D_{\mathbf{i}}$. In the random construction, the sets $D_{\mathbf{i}}$ are independent with a common, arbitrary distribution on $\mathcal A$.

For $j\in\mathcal I$, write $T_j=(D_{j\mathbf{i}})_{\mathbf{i}\in\mathcal I^*}$. Then
\[
    \nu_T=\frac1{|D_\varnothing|}\sum_{j\in D_\varnothing}
        \tau_{j/M}S_{1/M}\nu_{T_j}.
\]
On $X=\mathcal T\times\mathcal I$, define
\[
    P_{(T,i)}=\frac1{|D_\varnothing|}\sum_{j\in D_\varnothing}\delta_{(T_j,j)},
    \qquad \pi(T,i)=i,
    \qquad \lambda(T,i)=M^{-1}.
\]
The space $X$ is compact and $P$ is continuous. By \eqref{self_similarity_of_Markov_self_similar_measures}, $\mu_{(T,i)}=\tau_i\nu_T$, so $\mu_{(T,0)}=\nu_T$. Hence, the random measure $\nu_T$ is embedded into the family of Markov self-similar measures generated by $(X,P,\pi,\lambda)$.

Fractal percolation with retention probability $p\in(M^{-d},1)$, conditioned on survival, fits this construction after discarding cubes with no infinite line of descendants; see \cite[Proposition 5.28(i)]{LP16}. The associated random measure splits mass equally among the children with infinite descendants.

Another class of random fractal measures studied in the literature are random perturbations of self-similar measures, where the translations or the contractions in the IFS are replaced by random variables; see for instance \cite{PSS06}.

These are just two classes of examples of random fractal measures that can be represented as Markov self-similar measures, but hopefully they illustrate the richness and generality of the Markov self-similar measure framework.

\subsection{Furstenberg CP-processes}

Furstenberg's dyadic CP-processes \cite{Fur08} also fit naturally into the Markov self-similar framework. Given $Q\in\D_1([0,1)^d)$, let $h_Q$ be the homothety from $Q$ onto $[0,1)^d$. For a probability measure $\nu$ on $[0,1]^d$ and a dyadic cube $Q$ with $\nu(Q)>0$, define the magnified conditional measure
\[
M_Q\nu=h_Q\left(\frac{\nu|_Q}{\nu(Q)}\right)\in \Prob([0,1]^d).
\]
Let $\mathcal Y\subset \Prob([0,1]^d)$ be a Borel family (in the weak$^*$ topology) of probability measures on $[0,1]^d$ such that
\begin{enumerate}
\renewcommand{\labelenumi}{\rm{(\roman{enumi})}}
\item $\nu(\partial Q)=0$ for every $\nu\in \mathcal Y$ and every $Q\in\D_1([0,1)^d)$,
\item if $\nu\in\mathcal Y$ and $\nu(Q)>0$, then $M_Q\nu\in\mathcal Y$.
\end{enumerate}

We define a Markov self-similar model as follows. Let
\[
    X=\{(\nu,Q): \nu\in\mathcal Y,\ Q\in\D_1([0,1)^d),\ \nu(Q)>0\}.
\]
Define a Markov kernel on $X$ by
\[
    P_{(\nu,Q)}
    =
    \sum_{Q'\in\D_1([0,1)^d)} M_Q\nu(Q')\,\delta_{(M_Q\nu,Q')}.
\]
In other words, from the state $(\nu,Q)$ one passes to the magnified measure $M_Q\nu$ and then chooses the next dyadic child according to that measure. Writing $Q=a_Q+[0,1/2)^d$, define
\[
    \lambda(\nu,Q)=\frac{1}{2},
    \qquad
    \pi(\nu,Q)=a_Q.
\]
Then the affine map $\tau_{\pi(\nu,Q)}S_{\lambda(\nu,Q)}$ is exactly the inverse of $h_Q$. For $(\nu,Q)\in X$, set
\[
    \mu_{(\nu,Q)}=\tau_{\pi(\nu,Q)}S_{\lambda(\nu,Q)}M_Q\nu=\frac{\nu|_Q}{\nu(Q)}.
\]
Since
\[
    M_Q\nu
    =
    \sum_{Q'\in\D_1([0,1)^d)} M_Q\nu(Q')\, \tau_{\pi(M_Q\nu,Q')}S_{1/2} M_{Q'}(M_Q\nu),
\]
it follows after applying $\tau_{\pi(\nu,Q)}S_{1/2}$ that
\[
    \mu_{(\nu,Q)}
    =
    \sum_{Q'\in\D_1([0,1)^d)} M_Q\nu(Q')\,
    \tau_{\pi(\nu,Q)}S_{1/2}\mu_{(M_Q\nu,Q')}.
\]
This tells us that $(\mu_{(\mu,Q)})_{(\nu,Q)\in X}$ is the family of Markov self-similar measures generated by the above model $(X,P,\pi,\lambda)$.
In particular, every $\nu\in\mathcal Y$ decomposes as
\[
    \nu=\sum_{Q\in\D_1([0,1)^d)}\nu(Q)\mu_{(\nu,Q)}.
\]
Thus dyadic CP-processes provide a broad class of homogeneous Markov self-similar models, with the state encoding the current scenery and the chosen dyadic cell.

Note that in Furstenberg's work \cite{Fur08}, a central role is played by certain stationary measures of the Markov chain on $X$ - these are properly the CP-processes, rather than just the Markov dynamics itself.


\section{\texorpdfstring{$L^q$}{Lq} dimension of Markov self-similar measures}

\label{sec:Lq-dimension-Markov-SSM}

\subsection{The worst case spectrum \texorpdfstring{$T(q)$}{Tq}}

\label{S:worst-case-spectrum}

Fix a Markov self-similar model $\mathcal{X}=(X,P,\pi,\lambda)$. We write the estimates for $\supp(\mu_{x,n})\subset[0,1)^d$; they extend to uniformly bounded supports as explained in \S\ref{S:support-renormalization}.

Fix $q\in (1,\infty)$. Our goal is to understand the $L^q$ dimension of $\mu_x$ for $q>1$, under suitable assumptions. We start by recalling the relevant definitions, introducing notation along the way.

For a nonzero finite Borel measure $\nu$ on $\R^d$ with the bounded support, its \emph{level-$m$ dyadic discretization} is the atomic measure
\[
    \nu^{(m)}=\sum_{k\in\Z^d}\nu\bigl(2^{-m}(k+[0,1)^d)\bigr)\delta_{2^{-m}k},
\]
obtained by placing the mass of each cube at its corner with minimal coordinates. For an atomic measure $\rho$ (for example, $\rho=\nu^{(m)}$), we write $\|\rho\|_q=(\sum_{x\in\supp\rho}\rho(x)^q)^{1/q}$. With this notation, the $L^q$ spectrum of $\nu$ is given by
\[
\tau_\nu(q)=-\lim_{m\to\infty} \frac{\log \|\nu^{(m)}\|_q^q}{m}=-\lim_{m\to\infty}\frac{1}{m}\log\sum_{I\in\D_m}\nu(I)^q,
\]
provided the limit exists. In general, we define the upper and lower $L^q$ spectra by
\[
\overline{\tau}_\nu(q)=\limsup_{m\to\infty}\left(-\frac{\log \|\nu^{(m)}\|_q^q}{m}\right)
\quad\text{and}\quad
\underline{\tau}_\nu(q)=\liminf_{m\to\infty}\left(-\frac{\log \|\nu^{(m)}\|_q^q}{m}\right).
\]
Whenever the $L^q$ spectrum exists, we define the $L^q$ dimension of $\nu$ by
\[
    D_q(\nu)=\frac{\tau_\nu(q)}{q-1}.
\]
This normalization ensures that $D_q(\nu)\in [0,d]$ for every $q>1$.

For $x\in X$ and $m\in\N$, let
\[
\phi_m(x)=\|\mu_{x}^{(m)}\|_q^q.
\]
The map $x\mapsto\mu_x$ is a Borel kernel, so each $\phi_m$ is a positive Borel function.
We define the \emph{worst case $L^q$ spectrum} of the model $\mathcal{X}$ by
\begin{equation}\label{def:Tq}
T(q)=T_{\mathcal{X}}(q)=-\lim_{m\to\infty}\frac{1}{m}\,\sup_{y\in X}\log\phi_m(y).
\end{equation}
We will see in Lemma \ref{lem:Tq} below that the limit in \eqref{def:Tq} exists, and we will derive some additional properties of $T(q)$. For now, the main observation is that, by definition, 
\[
    \overline{\tau}_{\mu_x}(q)\ge \underline{\tau}_{\mu_x}(q)\ge T(q)
\qquad\text{for every }x\in X.
\]
In fact, this holds uniformly in $x$ in the sense of the $\liminf$ in the definition of $\underline{\tau}_{\mu_x}(q)$.

In many cases, $L^q$ norms are known to satisfy a sub-multiplicative relation along scales. For self-conformal measures, this was shown by Peres and Solomyak \cite{PS00}. In the dynamically-driven case of \cite{Shm19,CS24}, the $L^q$ norms satisfy a sub-multiplicative relation of the form
\begin{equation} \label{eq:submultiplicativity-dynamically-driven}
\phi_{m+n}(y)\lesssim \phi_m(y)\phi_n(T^{\lfloor m/\log(1/\lambda)\rfloor} y).
\end{equation}
This is a fundamental fact in \cite{Shm19,CS24} that, together with unique ergodicity of the underlying dynamics, allows to conclude that the ``worst case'' $L^q$ spectrum $T(q)$ is actually the $L^q$ spectrum of $\mu_y$ for almost every $y\in X$ (with respect to the unique invariant measure of $T$).

To our knowledge, there is no suitable analog of \eqref{eq:submultiplicativity-dynamically-driven} for general Markov self-similar measures, even under strong assumptions on the Markov chain like continuity, homogeneity and unique ergodicity. However, Corollary \ref{cor:subadditivity-with-max-over-orbit} yields a very weak substitute that, as we will see, is sufficient in many applications, and in particular recovers the main applications of \cite{Shm19,CS24} without any appeal to \eqref{eq:submultiplicativity-dynamically-driven}.

\subsection{Technical preliminaries}

We recall some basic estimates on $L^q$ norms for future reference.
\begin{lem}[\cite{Shm19}, Lemmas 4.1--4.3]\label{lem:partition-comparison}
Let $q\in(1,\infty)$.
\begin{enumerate}[label=\textup{(\roman*)}]
\item Let $\nu$ be a finite measure on a measurable space $Y$, and let $\mathcal P,\mathcal Q$ be countable families of measurable subsets of $Y$. Suppose that, for some $M\in\N$, each member of $\mathcal P$ is covered by at most $M$ members of $\mathcal Q$, and each member of $\mathcal Q$ intersects at most $M$ members of $\mathcal P$. Then
\[
    \sum_{E\in\mathcal P}\nu(E)^q\le M^q\sum_{F\in\mathcal Q}\nu(F)^q.
\]
In particular, if $\mathcal P,\mathcal Q$ are partitions of $Y$ and each member of either partition intersects at most $M$ members of the other, then
\[
    M^{-q}\sum_{F\in\mathcal Q}\nu(F)^q
    \le \sum_{E\in\mathcal P}\nu(E)^q
    \le M^q\sum_{F\in\mathcal Q}\nu(F)^q.
\]
\item Let $\nu=\sum_i\nu_i$ be a finite atomic measure. If each atom has positive mass for at most $M\in\N$ of the measures $\nu_i$, then
\[
    \sum_i\|\nu_i\|_q^q\le \|\nu\|_q^q\le M^{q-1}\sum_i\|\nu_i\|_q^q.
\]
\item For finite Borel measures $\nu,\eta$ on $\R^d$ and $m\in\N$,
\[
    \|(\nu*\eta)^{(m)}\|_q^q\sim_{d,q}\|\nu^{(m)}*\eta^{(m)}\|_q^q.
\]
\end{enumerate}
\end{lem}

\begin{lem}\label{lem:stopped-dyadic-comparison-stopped-section}
For every $m\in\N$ and $t\ge\max\{0,m-\alpha_{\max}\}$, we have $|\pi_\infty-\pi_{T_t}|_\infty\le2^{-t}\le\lambda_{\min}^{-1}2^{-m}$ uniformly on $X^{\N_0}$. Consequently, for every $x\in X$,
\[
\|\mu_x^{(m)}\|_q^q \sim_q \|\mu_{x,[t]}^{(m)}\|_q^q.
\]
\end{lem}
\begin{proof}
The bound follows from the series defining $\pi_\infty$. Each member of either pullback of $\D_m$ under $\pi_\infty$ and $\pi_{T_t}$ therefore intersects $O_{\lambda_{\min},d}(1)$ members of the other, up to a $P_x^\infty$-null set. The conclusion follows from Lemma \ref{lem:partition-comparison}(i).
\end{proof}

\begin{cor} \label{cor:subadditivity-with-max-over-orbit}
    For all $m,n\in\N$ and $x\in X$,
    \[
    \phi_{m+n}(x) \lesssim_q \phi_m(x) \|\phi_n\|_{L^\infty(Q_x^{[m]})}.
    \]
\end{cor}
\begin{proof}
By the Markov property at $T_m$ (obtained by taking $t\to\infty$ in (\ref{eq:stopped-self-similarity})),
\[
\mu_x=\int_{X^{\N_0}}\tau_{\pi_{T_m-1}}S_{\lambda(x_0^{T_m})}\mu_{x_{T_m}}\,dP_x^\infty.
\]
Fix $I\in\D_m$. Every component charging $I$ has $\pi_{T_m}\in3I$, since $\lambda_{\min} 2^{-m}\le\lambda(x_0^{T_m})<2^{-m}$. Jensen's inequality and Lemma \ref{lem:partition-comparison}(i) yield
\[
\sum_{J\in\D_{m+n}(I)}\mu_x(J)^q
\lesssim_q\mu_{x,[m]}(3I)^q\|\phi_n\|_{L^\infty(Q_x^{[m]})}.
\]
Summing over $I\in\D_m$ and applying Lemmas \ref{lem:partition-comparison}(i) and \ref{lem:stopped-dyadic-comparison-stopped-section}, we obtain
\[
\phi_{m+n}(x)\lesssim_q
\|\phi_n\|_{L^\infty(Q_x^{[m]})}
\sum_{I\in\D_m}\mu_{x,[m]}(3I)^q
\lesssim_q\phi_m(x)\|\phi_n\|_{L^\infty(Q_x^{[m]})}.
\qedhere
\]
\end{proof}

Next, we confirm that the limit in the definition of worst-case spectrum $T(q)$ exists, and we establish some of its basic properties.
\begin{lem} \label{lem:Tq}
The limit
\[
T(q)=-\lim_{m\to\infty} \frac{1}{m}\,\sup_{y\in X} \log\phi_m(y)
\]
exists. 
There exists a constant $C_q>0$ such that
\begin{equation}\label{eq:Tq-variational}
T(q)=\sup_{m\ge1}-\frac{\sup_{y\in X}\log\phi_m(y)+C_q}{m}.
\end{equation}
 The function $T(q)$ is concave and non-decreasing in $q$. Moreover, for every $\e>0$, if $m\in\N$ is sufficiently large, then
\[
\|\mu_x^{(m)}\|_q^q\le 2^{-m(T(q)-\e)}\quad\text{for every }x\in X.
\]
\end{lem}
\begin{proof}
Write $a_m(q)=\sup_{y\in X}\log\phi_m(y)$. By Corollary \ref{cor:subadditivity-with-max-over-orbit},
\begin{equation}\label{eq:an-almost-subadditive}
    a_{m+n}(q)\le a_m(q)+a_n(q)+C_q.
\end{equation}
Fekete's lemma applied to $a_m(q)+C_q$ establishes the limit and \eqref{eq:Tq-variational}. By the definition and the existence of the limit for $T(q)$, $\|\mu_x^{(m)}\|_q^q\le2^{-m(T(q)-\e)}$ for large $m$, uniformly in $x$. Each function $q\mapsto\log\phi_m(y)=\log\sum_{I\in\D_m}\mu_y(I)^q$ is convex and non-increasing, as is their supremum $a_m(q)$. Passing to the limit shows that $T(q)$ is concave and non-decreasing.
\end{proof}

\subsection{Large \texorpdfstring{$L^q$}{Lq} norms along orbits}

For $m\in\N$, a \emph{stopped scale-$m$ orbit} is a sequence $(y_0,\ldots,y_{\ell-1})$ of the form $y_j=x_{t_j}$ along a Markov orbit $(x_n)_{n\ge0}$, where $t_0=0$ and
\[
t_{j+1}=\min\{n>t_j:\lambda(x_{t_j}^n)<2^{-m}\}
\qquad(0\le j<\ell-1).
\]

We iterate Corollary \ref{cor:subadditivity-with-max-over-orbit} by conditioning successive stopped path laws.

\begin{cor} \label{cor:iterating-subadditivity}
    For all $m,\ell\in\N$ and $y_0\in X$, there is a stopped scale-$m$ orbit $(y_0,\ldots,y_{\ell-1})$ with $y_{j+1}\in \supp(Q_{y_j}^{[m]})$ for $0\le j<\ell-1$, such that
    \[
     \phi_{\ell m}(y_0) \le C_q^\ell\prod_{j=0}^{\ell-1}\phi_m(y_j).
    \]
\end{cor}
\begin{proof}
Choose $C_q\ge1$ so that Corollary \ref{cor:subadditivity-with-max-over-orbit} holds with constant $C_q/2$. For $r=0,\ldots,\ell-2$, sample from $P_{y_r}^{[m]}$ conditioned on
\[
\phi_{(\ell-r-1)m}(x_{T_m})
\ge\frac{\phi_{(\ell-r)m}(y_r)}{C_q\phi_m(y_r)},
\]
and set $y_{r+1}=x_{T_m}$. Each conditioning event has positive probability by Corollary \ref{cor:subadditivity-with-max-over-orbit}. Concatenating these Borel conditional kernels and continuing with $P$ defines a law absolutely continuous with respect to $P_{y_0}^\infty$. Almost every trajectory is a Markov orbit with $y_{r+1}\in \supp(Q_{y_r}^{[m]})$. Multiplying the conditioning inequalities yields the claim.
\end{proof}

For a set $E\subset\{0,1,\ldots,N\}$, we define its \emph{Schnirelmann density} by
\[
\underline{d}_S(E)=\min_{1\le j\le N+1}\frac{|E\cap\{0,\ldots,j-1\}|}{j}.
\]
Note that this depends on both $E$ and $N$, but $N$ will always be clear from context.

We record a standard maximal-partial-sum argument, closely related to the cycle lemma. It will allow us to upgrade ``positive lower density'' to ``positive Schnirelmann density'' later on.
\begin{lem} \label{lem:combinatorial}
Fix $N\in\N$ and $0<\delta<\eta<1$. Let $E\subset\{0,1,\ldots,N\}$ be a set with cardinality $\ge(1-\delta)(N+1)$. Then there is some integer $0\le j\le\delta(N+1)/\eta$ such that the tail $(E\cap\{j,\ldots,N\})-j\subset\{0,1,\dots, N-j\}$ has Schnirelmann density at least $1-\eta$.
\end{lem}
\begin{proof}
Let $z_n=|\{0\le i<n:i\notin E\}|$ for $0\le n\le N+1$, and choose $j$ maximizing $z_j-\eta j$. Since $z_0=0$, we have $\eta j\le z_j\le z_{N+1}\le\delta(N+1)$. For every $1\le n\le N+1-j$, maximality yields $z_{j+n}-z_j\le\eta n$, so at least $(1-\eta)n$ of the integers $j,\ldots,j+n-1$ belong to $E$.
\end{proof}

\begin{prop} \label{prop:good-step-n-orbit}
    For every $\e,\eta>0$, all sufficiently large $m$, and every $k\in\N$, there is a stopped scale-$m$ orbit $(y_0,\ldots,y_{k-1})$ with $y_{i+1}\in\supp(Q_{y_i}^{[m]})$ for $0\le i<k-1$, such that
    \[
        \underline{d}_S\left\{ i\in \{0,\ldots,k-1\}: \phi_m(y_i) \ge 2^{-m(T(q)+\e)} \right\} \ge 1-\eta.
    \]
\end{prop}
\begin{proof}
    It suffices to consider $0<\eta<1/2$. 
Let $C_q$ be the constant in Corollary \ref{cor:iterating-subadditivity}. For sufficiently large $m$, Lemma \ref{lem:Tq} ensures that $\log C_q\le m\e\eta^2/16$, that $\phi_m(y)\le2^{-m(T(q)-\e\eta^2/16)}$ for all $y$, and that, for every $k\ge1$, there is $z_0$ with $\phi_{2km}(z_0)\ge2^{-2km(T(q)+\e\eta^2/16)}$. Applying Corollary \ref{cor:iterating-subadditivity}, we obtain a stopped scale-$m$ orbit $(z_0,\ldots,z_{2k-1})$ with $z_{j+1}\in\supp(Q_{z_j}^{[m]})$. Let 
\[
E=\{j:\phi_m(z_j)<2^{-m(T(q)+\e)}\}.
\] 
Then
\[
2^{-2km(T(q)+\e\eta^2/16)}
\le C_q^{2k}\prod_{j=0}^{2k-1}\phi_m(z_j)
\le2^{-2kmT(q)+m\e\eta^2k/4-m(\e+\e\eta^2/16)|E|}.
\]
Taking logarithms yields $|E|\le\eta^2k/2$.

    Apply Lemma \ref{lem:combinatorial}, with $\delta=\eta^2/4$, to the set $\{0,\ldots,2k-1\}\setminus E$. We obtain $j_0\le\eta k/2$ such that
    \[
        \bigl|\{0\le i<n:\phi_m(z_{j_0+i})\ge2^{-m(T(q)+\e)}\}\bigr|
        \ge(1-\eta)n
        \qquad(1\le n\le 2k-j_0).
    \]
    Since $2k-j_0\ge k$, the orbit $y_i=z_{j_0+i}$, $0\le i<k$, has the required property.
\end{proof}

\subsection{\texorpdfstring{$L^q$}{Lq}-flattening and \texorpdfstring{$L^q$}{Lq} norms at finer scales}
In this subsection we establish Proposition \ref{prop:Lq-at-finer-scale}, a finitary Markov analog of \cite[Proposition 5.2]{Shm19}. We assume the $L^q$ flattening condition from Definition \ref{dfn:Lq-flattening}.

Given a dyadic cube $I\in\D_m$, we write $h_I$ for the homothety that renormalizes $I$ to $[0,1)^d$. Fix $x\in X$, $m\in\N$, and $I\in\D_m$. When $\mu_{x,[m]}(I)>0$, let $\overline P_{x,I}^\infty$ denote the law of the full Markov orbit started from $x$, conditioned on the event $\pi_{T_m}\in I$. We define
\begin{align*}
\mu_{x,I}&=h_I\pi_{T_m}\overline P_{x,I}^\infty=h_I\left(\frac{\mu_{x,[m]}|_I}{\mu_{x,[m]}(I)}\right),\\
\nu_{x,I,[n]}&=h_I\int_{X^{\N_0}}\tau_{\pi_{T_m-1}(x_0^\infty)}S_{\lambda(x_0^{T_m})}\mu_{x_{T_m},[n]}\,d\overline P_{x,I}^\infty(x_0^\infty).
\end{align*}
If $\mu_{x,[m]}(I)=0$, we define $\mu_{x,I}$ and $\nu_{x,I,[n]}$ arbitrarily, since every subsequent estimate is weighted by $\mu_{x,[m]}(I)^q$.

\begin{lem} \label{lem:Lq-simple-comparison}
For every $\rho>0$ and all sufficiently large $n$,
\[
\|\nu_{x,I,[n]}^{(n)}\|_q^q\le 2^{-(T(q)-\rho)n}
\]
uniformly in $x\in X$, $m\in\N$, and $I\in\D_m$ with $\mu_{x,[m]}(I)>0$.
\end{lem}
\begin{proof}
For $\overline P_{x,I}^\infty$-almost every orbit, the inverse image of $\D_n$ under $h_I\tau_{\pi_{T_m-1}(x_0^\infty)}S_{\lambda(x_0^{T_m})}$ is a partition into cubes of side length between $2^{-n}$ and $\lambda_{\min}^{-1}2^{-n}$. Jensen's inequality and Lemmas \ref{lem:partition-comparison}(i) and \ref{lem:stopped-dyadic-comparison-stopped-section} yield
\[
\|\nu_{x,I,[n]}^{(n)}\|_q^q
\lesssim_q\int\|\mu_{x_{T_m},[n]}^{(n)}\|_q^q\,d\overline P_{x,I}^\infty
\lesssim_q\int\phi_n(x_{T_m})\,d\overline P_{x,I}^\infty.
\]
The bound $\|\nu_{x,I,[n]}^{(n)}\|_q^q\le2^{-n(T(q)-\rho)}$, uniform in $x,m,I$, follows from Lemma \ref{lem:Tq}.
\end{proof}

Heuristically, one may expect that when $\mu_{x,I}$ is ``exponentially non-atomic at scale $2^{-n}$'', the estimate from Lemma \ref{lem:Lq-simple-comparison} can be improved. This is because $\nu_{x,I,[n]}$ is obtained by averaging renormalized stopped descendants of scale about $2^{-n}$, and in many concrete settings this average has a convolution structure that flattens $L^q$ norms exponentially. This motivates the following key definition, which abstracts away the additive-combinatorial input needed to obtain $L^q$ dimension estimates.
\begin{dfn} \label{dfn:Lq-flattening}
  We say that the Markov self-similar model $(X,P,\pi,\lambda)$ satisfies the \emph{$L^q$ flattening condition} for $q>1$ if the following holds. For every $\sigma>0$ there are $\e=\e(\sigma)>0$ and $n_0=n_0(\sigma)\in\N$ such that, whenever $x\in X$, $I\in\D_m$, $n\ge m\geq n_0$, and $\mu_{x,[m]}(I)>0$,
  \[
    \|\mu_{x,I}^{(n)}\|_q^q \le 2^{-\sigma n}
    \Longrightarrow
    \|\nu_{x,I,[n]}^{(n)}\|_q^q \le 2^{-(T(q)+\e)n}.
  \]
  We refer to the function $\e(\sigma)$ as the \emph{$L^q$-flattening rate} of the model.
\end{dfn}

The following proposition is a Markov, finitary analog of \cite[Proposition 5.2]{Shm19}. It is the key consequence of the $L^q$ flattening condition.
\begin{prop} \label{prop:Lq-at-finer-scale}
    Suppose the Markov self-similar model satisfies the $L^q$ flattening condition for some $q>1$ with flattening rate $\e(\sigma)$. Fix an integer $R\ge 2$ and $\sigma>0$, and write $\e=\e(\sigma)\le\sigma$. Suppose $m$ is sufficiently large in terms of $q,R,\sigma$, and suppose $x\in X$ is such that
    \begin{equation} \label{eq:Lq-almost-sharp-assumption}
    \phi_{(R+1)m}(x)\ge 2^{-(R+1)m(T(q)+\e/2)}.
    \end{equation}
    Then,
    \[
     \left\|\mu_{x,[m]}^{((R+1)m)}\right\|_q^q  \gtrsim 2^{-m(T(q)+2 R\sigma)}.
    \]
\end{prop}
\begin{proof}
For $n\in\N$, the Markov property at $T_m$ shows that $\sum_{I\in\D_m}\mu_{x,[m]}(I)\,h_I^{-1}\nu_{x,I,[n]}=\int_{X^{\N_0}}\tau_{\pi_{T_m-1}(x_0^\infty)}S_{\lambda(x_0^{T_m})}\mu_{x_{T_m},[n]}\,dP^\infty_x$ is the law of $\pi_{T_m+T_n(x_{T_m}^\infty)}$ under $P_x^\infty$. Each summand is supported in $3I$, and 
\[
|\pi_\infty-\pi_{T_m+T_n(x_{T_m}^\infty)}|_\infty\le2^{-(m+n)}
\] 
uniformly on $X^{\N_0}$. Lemma \ref{lem:partition-comparison}(i)--(ii) therefore yields
\begin{equation}\label{eq:Lq-two-scale-cont}
\phi_{m+n}(x)\sim_q\sum_{I\in\D_m}\mu_{x,[m]}(I)^q\cdot\|\nu_{x,I,[n]}^{(n)}\|_q^q.
\end{equation}

Take $n=Rm$ and $m\geq n(\sigma)$, and
let $\D'=\{I\in\D_m:\|\mu_{x,I}^{(Rm)}\|_q^q\le2^{-\sigma Rm}\}$. The flattening condition bounds $\|\nu_{x,I,[Rm]}^{(Rm)}\|_q^q$ by $2^{-Rm(T(q)+\e)}$ on $\D'$. From \eqref{eq:Lq-two-scale-cont}, Lemmas \ref{lem:Lq-simple-comparison} and \ref{lem:Tq}, and $\sum_{I\in\D'}\mu_{x,[m]}(I)^q\lesssim_q\phi_m(x)$, we deduce that
\[
\phi_{(R+1)m}(x)\lesssim_q
2^{-m((R+1)T(q)+R\e-\e/8)}
+2^{-Rm(T(q)-\e/8)}\sum_{I\in\D_m\setminus\D'}\mu_{x,[m]}(I)^q.
\]
Since $R\ge2$, it follows from \eqref{eq:Lq-almost-sharp-assumption} that the contribution from $\D'$ is $o(\phi_{(R+1)m}(x))$. Hence
\[
\sum_{I\in\D_m\setminus\D'}\mu_{x,[m]}(I)^q
\gtrsim_q2^{Rm(T(q)-\e/8)}\phi_{(R+1)m}(x)
\ge2^{-m(T(q)+R\sigma)},
\]
using $\e\le\sigma$. By the definition of $\mu_{x,I}$,
\[
\|\mu_{x,[m]}^{((R+1)m)}\|_q^q=\sum_{I\in\D_m}\mu_{x,[m]}(I)^q\|\mu_{x,I}^{(Rm)}\|_q^q
\ge2^{-Rm\sigma}\sum_{I\in\D_m\setminus\D'}\mu_{x,[m]}(I)^q
\gtrsim_q2^{-m(T(q)+2R\sigma)}.
\qedhere
\]
\end{proof}

\section{On \texorpdfstring{$L^q$}{Lq} flattening}\label{S:Lq-flattening}

\subsection{Sub-exponential endpoint growth and \texorpdfstring{$L^q$}{Lq} flattening}\label{S:subexponential-orbit-growth}

In general, verifying $L^q$ flattening can be very challenging. We hope the $L^q$ porosity approach introduced by the second author \cite{Usu25} and developed further by H.~Ren and the authors \cite{RSU26} can be used to verify $L^q$ flattening in many cases. In this section, we consider Markov self-similar models with sub-exponential endpoint growth. While this condition is quite restrictive, it is satisfied in the examples that are relevant to our main results, and many others.

For $x\in X$ and $t\ge0$, let
\[
    \mathcal O_t(x)=\left\{\bigl(\lambda(x_0^n),x_n\bigr):
    x_0^n\in\supp(P_x^{[t]})\right\}.
\]
We remark that $\mathcal O_t(x)$ consists of pairs recording both the contraction ratio and the stopped state.
\begin{dfn}\label{dfn:subexponential-orbit-growth}
We say that the Markov self-similar model has \emph{sub-exponential endpoint growth} if $\sup_{x\in X}|\mathcal O_m(x)|<\infty$ for every $m\in\N$ and
\[
    \lim_{m\to\infty}\frac{1}{m}\sup_{x\in X}\log|\mathcal O_m(x)|=0.
\]
\end{dfn}

The next lemma shows that, under sub-exponential endpoint growth, $L^q$ flattening reduces to checking flattening for convolutions with a single measure.
\begin{lem}\label{lem:subexponential-orbit-growth-flattening}
Fix $q>1$, and let $\mathcal X=(X,P,\pi,\lambda)$ be a Markov self-similar model on $\R^d$ with sub-exponential endpoint growth. Suppose that, for every $\sigma>0$ and $C\ge1$, there are $\eta>0$ and $n_0\in\N$ such that
\begin{equation}\label{eq:ordinary-tail-flattening}
    \bigl\|(\theta*S_r\mu_{x,[n]})^{(n)}\bigr\|_q^q
    \le 2^{-(T(q)+\eta)n}
\end{equation}
for every $n\ge n_0$, $x\in X$, $r\in[C^{-1},C]$, and every probability measure $\theta$ supported on $[-C,C]^d$ with $\|\theta^{(n)}\|_q^q\le2^{-\sigma n}$. Then $\mathcal X$ satisfies the $L^q$ flattening condition.
\end{lem}
\begin{proof}
Fix $\sigma>0$, $x\in X$, $I\in\D_m$ with $\mu_{x,[m]}(I)>0$, and $n\ge m$, and assume $\|\mu_{x,I}^{(n)}\|_q^q\le2^{-\sigma n}$. Partition the paths under $\overline P_{x,I}^\infty$ by $(2^m\lambda(x_0^{T_m}),x_{T_m})=(r,z)\ (r\in[\lambda_{\min},1), z\in X)$. Let $p_{r,z}$ be the probability of each class and $\theta_{r,z}$ the conditional law of $h_I\pi_{T_m-1}$, omitting classes of zero probability. There are at most $|\mathcal O_m(x)|$ classes. By the definitions of $\mu_{x,I}$ and $\nu_{x,I,[n]}$,
\begin{equation}\label{eq:orbit-growth-class-decompositions}
    \mu_{x,I}=\sum_{r,z}p_{r,z}\tau_{r\pi(z)}\theta_{r,z},
    \qquad
    \nu_{x,I,[n]}=\sum_{r,z}p_{r,z}\theta_{r,z}*S_r\mu_{z,[n]}.
\end{equation}
Lemma \ref{lem:partition-comparison}(i) and the basic fact $(s+t)^q\geq s^q+t^q$ yield
\begin{equation}\label{eq:orbit-growth-class-Lq-sum}
    \sum_{r,z}p_{r,z}^q\|\theta_{r,z}^{(n)}\|_q^q
    \lesssim_q\|\mu_{x,I}^{(n)}\|_q^q\le2^{-\sigma n}.
\end{equation}

Since $r\in[\lambda_{\min},1)$ and $\theta_{r,z}$ is supported on $[0,1]^d-r\pi(z)$, we may apply \eqref{eq:ordinary-tail-flattening} with $\sigma/2$ and $C=\max\{\lambda_{\min}^{-1},1+\sup_{x\in X}|\pi(x)|_\infty\}$. This provides some $\eta>0$ and $n_0\in\N$ such that
\[
    \|(\theta_{r,z}*S_r\mu_{z,[n]})^{(n)}\|_q^q
    \le2^{-(T(q)+\eta)n}
    \quad\text{if }\|\theta_{r,z}^{(n)}\|_q^q\le2^{-\sigma n/2},\quad\text{assuming }n\geq n_0.
\]
For the remaining pairs, we deduce from H\"older's inequality and \eqref{eq:orbit-growth-class-Lq-sum} that
\[
    \sum_{\substack{r,z\\\|\theta_{r,z}^{(n)}\|_q^q>2^{-\sigma n/2}}}p_{r,z}
    \lesssim_q|\mathcal O_m(x)|^{1-1/q}2^{-\sigma n/(2q)}
    \le2^{-\sigma n/(4q)}
\]
for sufficiently large $m$ in terms of $q$ and $\sigma$, uniformly in $m\le n$ and $x\in X$, by Definition \ref{dfn:subexponential-orbit-growth}. For every pair, Young's inequality and Lemmas \ref{lem:partition-comparison}(i), (iii), \ref{lem:stopped-dyadic-comparison-stopped-section}, and \ref{lem:Tq} yield
\[
    \|(\theta_{r,z}*S_r\mu_{z,[n]})^{(n)}\|_q^q
    \le2^{-(T(q)-\sigma/(8q))n}
\]
for sufficiently large $n$. Jensen's inequality in \eqref{eq:orbit-growth-class-decompositions} therefore yields
\[
    \|\nu_{x,I,[n]}^{(n)}\|_q^q
    \le2^{-(T(q)+\eta)n}
       +2^{-(T(q)+\sigma/(8q))n}
    \le2^{-(T(q)+\frac12\min\{\eta,\sigma/(8q)\})n}
\]
for sufficiently large $n\ge m$, uniformly in $x$ and $I\in\D_m$ with $\|\mu_{x,I}^{(n)}\|_q^q\leq2^{-\sigma n}$, as required.
\end{proof}

\subsection{Flattening of Markov self-similar measures under convolution}

Recall that $\mathbb{G}(d,d-1)$ is the Grassmannian of hyperplanes in $\R^d$. Following \cite{CS24}, we say that a Markov self-similar model $\mathcal{X}$ on $\R^d$ is \emph{$q$-unsaturated on lines} if
\begin{equation}\label{eq:q-unsaturation}
    T_{\mathcal X}(q)<T_{\Pi_V\mathcal X}(q)+q-1 \quad\text{for all } V\in\mathbb{G}(d,d-1).
\end{equation}
Note that if \eqref{eq:q-unsaturation} holds, then $T_{\mathcal X}(q)<d(q-1)$, since $T_{\Pi_V\mathcal X}(q)\le(d-1)(q-1)$ for all $V$. In dimension $d=1$, unsaturation on lines is equivalent to $T_{\mathcal X}(q)<q-1$.

As in \cite{CS24}, this pointwise definition has a uniform version:
\begin{lem} \label{lem:uniform-unsaturation}
    The map $V\mapsto T_{\Pi_V\mathcal X}(q)$ is lower semicontinuous on $\mathbb{G}(d,d-1)$. Consequently, if $\mathcal X$ is $q$-unsaturated on lines, then $\gamma_{\mathcal{X}}(q):=\inf_{V\in\mathbb{G}(d,d-1)}T_{\Pi_V\mathcal X}(q)$ satisfies $\gamma_{\mathcal{X}}(q)>T_{\mathcal X}(q)-(q-1)$, and for every $\rho>0$ and all sufficiently large $n$,
\[
    \|(\Pi_W\mu_x)^{(n)}\|_q^q
    \le 2^{-(\gamma_{\mathcal{X}}(q)-\rho)n}
    \quad\text{for all } W\in\mathbb{G}(d,d-1),\ x\in X.
\]
\end{lem}
\begin{proof}
Let
\[
    a_m(V)=\sup_{x\in X}\log\|(\Pi_V\mu_x)^{(m)}\|_q^q,
\]
using orthonormal coordinates on $V$. Corollary \ref{cor:subadditivity-with-max-over-orbit} and \eqref{eq:Tq-variational} yield
\[
    a_{m+n}(V)\le a_m(V)+a_n(V)+C_q,
    \qquad
    T_{\Pi_V\mathcal X}(q)=\sup_{m\ge1}-\frac{a_m(V)+C_q}{m},
\]
where the constant $C_q$ is independent of $V$. After enlarging $C_q$, Lemma \ref{lem:partition-comparison}(i) also shows that  $a_m(W)\le a_m(V)+C_q$ whenever $\|\Pi_W-\Pi_V\|\le 2^{-m}$. Letting $m\to\infty$, we deduce the claimed semicontinuity. The strict inequality for $\gamma_{\mathcal X}(q)$ follows by compactness of $\mathbb{G}(d,d-1)$, and the bound for $\|(\Pi_V\mu_x)^{(n)}\|_q^q$ follows from Corollary \ref{cor:subadditivity-with-max-over-orbit} and Lemma \ref{lem:partition-comparison}(i) (combined with the idea of the proof of Fekete's lemma), and compactness of $\mathbb{G}(d,d-1)$.
\end{proof}

\begin{thm} \label{thm:markov-convolution-flattening}
        Let $\mathcal X$ be a Markov self-similar model on $\R^d$ and fix $q>1$ such that $\mathcal X$ is $q$-unsaturated on lines. 
   
        For every $\sigma>0$ and $C\ge1$, there are $\e>0$ and $M_0\in\N$ such that, for every $M\ge M_0$, $x\in X$, $r\in[C^{-1},C]$, and probability measure $\theta$ supported on $[-C,C]^d$ with $\|\theta^{(M)}\|_q^q\le2^{-\sigma M}$,
    \[
        \bigl\|(\theta*S_r\mu_{x,[M]})^{(M)}\bigr\|_q^q
        \le 2^{-(T_{\mathcal X}(q)+\e)M}.
    \]
\end{thm}

For the proof we need the higher-dimensional inverse theorem of the first author \cite[Theorem 1.2]{Shm25}, in the form of \cite[Theorem 2.1]{CS24}.
\begin{thm}\label{thm:inverse-thm}
Fix $d\in\N$, $q>1$ and $\delta>0$. There is $L_0=L_0(d,q,\delta)\in\N$ such that, for every $L\ge L_0$, there are $\e_0>0$ and $S_0\in\N$ with the following property. Let $S\ge S_0$, $M=SL$, and let $\nu,\theta$ be probability measures supported on $2^{-M}\Z^d\cap[1/3,2/3)^d$. If
\[
    \|\nu*\theta\|_q\ge2^{-\e_0M}\|\nu\|_q,
\]
then, after translating $\nu$ by elements of $2^{-M}\Z^d\cap[-1/3,1/3]^d$, there are $F\subset\supp(\nu)$ and $\mathcal S_1\subset\{0,\ldots,S-1\}$ such that:
\begin{enumerate}[label=\textup{(\roman*)}]
\item $\|\nu|_F\|_q\ge2^{-\delta M}\|\nu\|_q$, and $\nu(a)\le2\nu(b)$ for every $a,b\in F$.
\item For each $s<S$, the number $R_s=|\D_{(s+1)L}(F\cap Q)|$ is independent of $Q\in\D_{sL}(F)$, and $F\cap Q\subset\frac12 Q$ for every such $Q$.
\item We have
\[
    |\mathcal S_1|L\ge\frac{-\log\|\theta\|_q^q-q\delta M-q}{d(q-1)}.
\]
For every $s\in\mathcal S_1$ and $Q\in\D_{sL}(F)$, there is $V\in\mathbb{G}(d,d-1)$ such that
\[
    |\D_{(s+1)L}(\Pi_V(F\cap Q))|
    \lesssim 2^{(\delta-1)L}|\D_{(s+1)L}(F\cap Q)|,
\]
with an implicit constant depending on $d$ only.
\end{enumerate}
\end{thm}

\begin{proof}[Proof of Theorem \ref{thm:markov-convolution-flattening}]
We adapt the refinement argument of \cite[\S3.3]{CS24}. The bound \eqref{eq:markov-convolution-flattening-limit} is vacuous for large $M$ if $\sigma>d(q-1)$, so we assume that $\sigma\le d(q-1)$.

It is enough to establish 
\begin{equation}\label{eq:markov-convolution-flattening-limit}
    \|\theta^{(M)}*\mu_x^{(M)}\|_q^q\le2^{-(T_{\mathcal X}(q)+\varepsilon)M}
\end{equation}
under the assumption that both measures are supported on $[1/3,2/3)^d$. Indeed, by several applications of Lemma \ref{lem:partition-comparison},  translating and rescaling the measures incurs only a constant factor in the $L^q$ norm, which can be absorbed into $\varepsilon$. Moreover, another application of Lemma \ref{lem:partition-comparison} together with Lemma \ref{lem:stopped-dyadic-comparison-stopped-section} show that 
\[
 \|\theta^{(M)}*\mu_x^{(M)}\|_q^q \sim \bigl\|\bigl(\theta*\mu_{x,[M]}\bigr)^{(M)}\bigr\|_q^q.
\]

In light of Lemma \ref{lem:uniform-unsaturation}, we may choose $\eta>0$ such that
\begin{equation}\label{eq:markov-uniform-unsaturation}
    4q\eta<\gamma_{\mathcal X}(q)+q-1-T_{\mathcal X}(q),
\end{equation}
and, for all sufficiently large $n$,
\begin{equation}\label{eq:markov-uniform-projected-Lq}
    \|(\Pi_V\mu_x)^{(n)}\|_q^q
    \le2^{-(\gamma_{\mathcal X}(q)-q\eta)n}\quad\text{for all }V\in\mathbb{G}(d,d-1), x \in X.
\end{equation}

Choose $\delta>0$ with 
\begin{equation}\label{eq:markov-delta-choice}
\delta<\min\{\eta/2,\sigma/(4q),\eta\sigma/(8dq(q-1))\}.
\end{equation}
Choose $L\ge L_0(q,\delta)$ from Theorem \ref{thm:inverse-thm} sufficiently large that \eqref{eq:markov-uniform-projected-Lq} holds at $n=L$, and such that
\begin{equation}\label{eq:markov-L-choice}
    L^{-1}<\eta\sigma/(8d(q-1)).
\end{equation}
With $L$ fixed, let $\e_0$ be the parameter from Theorem \ref{thm:inverse-thm}, and choose 
\begin{equation}\label{eq:markov-epsilon-choice}
0<\varepsilon<\min\{q\e_0/4,\eta\sigma/(8d(q-1))\}.
\end{equation}

Suppose $M=SL$ is large and \eqref{eq:markov-convolution-flattening-limit} fails. By Lemma \ref{lem:Tq}, $\|\mu_x^{(M)}\|_q^q\le2^{-(T_{\mathcal X}(q)-q\e_0/4)M}$, and so
\[
\frac{\|\theta^{(M)}*\mu_x^{(M)}\|_q}{\|\mu_x^{(M)}\|_q}
>2^{-(\varepsilon/q+\e_0/4)M}\stackrel{\eqref{eq:markov-epsilon-choice}}{\ge}2^{-\e_0M}.
\]
Apply Theorem \ref{thm:inverse-thm} to $\mu_x^{(M)}$ and $\theta^{(M)}$. Write $\widetilde\mu=\tau_t\mu_x$ for the translate of $\mu_x$ given by the theorem and let $F,(R_s)_{s<S}$ be the resulting set and branching numbers. Let $\mathcal S_1$ be the set of scales from Theorem \ref{thm:inverse-thm}(iii). We have
\begin{equation}\label{eq:markov-many-saturated-blocks}
    |\mathcal S_1|L\ge\frac{-\log\|\theta^{(M)}\|_q^q-q\delta M-q}{d(q-1)}
    \stackrel{\eqref{eq:markov-delta-choice}}{\ge}\frac{\sigma M}{2d(q-1)}.
\end{equation}
Construct $F'\subset F$ through the $S$ blocks, starting with $s=0$ and $F'=F$. For $s\notin\mathcal S_1$, leave $F'$ unchanged. For $s\in\mathcal S_1$ and each $Q\in\D_{sL}(F')$, let $\mathcal E=\D_{(s+1)L}(F\cap Q)$, so $|\mathcal E|=R_s$. Every member of $\mathcal E$ meets $\frac12 Q$ by part (ii). Since $\delta<\eta/2$ by \eqref{eq:markov-delta-choice}, Theorem \ref{thm:inverse-thm}(iii) gives a hyperplane $V$ such that, for large $L$,
\begin{equation}\label{eq:markov-local-line-saturation}
    \left|\D_{(s+1)L}(\Pi_V(\cup \mathcal E))\right|
    \stackrel{\eqref{eq:markov-delta-choice}}{\le}2^{-(1-\eta)L}|\mathcal E|.
\end{equation}
By the Markov property at $T_{sL+3}$,
\begin{equation}\label{eq:markov-stopped-decomposition}
    \widetilde{\mu}=\tau_t\mu_x=\int_{X^{\N_0}} \tau_{t+\pi_{T_{sL+3}-1}(x_0^\infty)}S_{\lambda(x_0^{T_{sL+3}})}\mu_{x_{T_{sL+3}}}\,dP_x^\infty.
\end{equation}
Every component charging $\bigcup \mathcal E$ is supported in $Q$, since its diameter in $|\cdot|_\infty$ is less than $2^{-(sL+3)}$ and each $J\in\mathcal E$ meets $\frac12 Q$. Let
$\Omega\subset X^{\N_0}$ be the set of $x_0^\infty$ such that the corresponding component measure charges $\bigcup\mathcal{E}$, which satisfies $P^\infty_x(\Omega)\leq\widetilde{\mu}(Q)$. Let
$\nu=\widetilde{\mu}|_{\bigcup \mathcal{E}}$. Then, apply Lemma \ref{lem:partition-comparison}(i) and \eqref{eq:markov-uniform-projected-Lq} to the projected component measures, then Jensen's inequality, to obtain
\begin{align}\label{eq:markov-local-projected-bound}
    \|(\Pi_V\nu)^{((s+1)L)}\|_q^q
    &\leq P^\infty_x(\Omega)^q\int_\Omega\|(\tau_{\Pi_V(t+\pi_{T_{sL+3}-1}(x_0^\infty))}S_{\lambda(x_0^{T_{sL+3}})}\Pi_V\mu_{x_{T_{sL+3}}})^{((s+1)L)}\|_q^q\ d(P^\infty_x)_\Omega\nonumber\\
    &\le 2^{-(\gamma_{\mathcal X}(q)-2q\eta)L}\widetilde{\mu}(Q)^q.
\end{align}
The constants are uniform in $V$ by \eqref{eq:markov-uniform-projected-Lq}.

Remove from $\mathcal E$ the $\lceil2^{-\eta L}|\mathcal E|\rceil$ cubes of largest $\nu$-mass, obtaining $\mathcal E'$. Then $\nu(J)\le2^{\eta L}\nu(\R^d)/|\mathcal E|$ for $J\in\mathcal E'$. H\"older's inequality on the projected cubes meeting $\Pi_V(\bigcup \mathcal E)$, followed by \eqref{eq:markov-local-line-saturation} and \eqref{eq:markov-local-projected-bound}, yields
\begin{align}
    \sum_{J\in\mathcal{E}'}\widetilde{\mu}(J)^q=\sum_{J\in\mathcal E'}\nu(J)^q
    &\le \left(\frac{2^{\eta L}}{|\mathcal E|}\right)^{q-1}\nu(\R^d)^q\notag\\
    &\le 2^{-(q-1)(1-2\eta)L}\|(\Pi_V\nu)^{((s+1)L)}\|_q^q\notag\\
    &\le 2^{-(T_{\mathcal X}(q)+\eta)L}\widetilde{\mu}(Q)^q.\label{eq:markov-local-unsaturation-Lq}
\end{align}
Replace $F'$ by $F'\setminus\bigcup(\mathcal E\setminus\mathcal E')$. Since $|\mathcal E\setminus\mathcal E'|\le \tfrac{1}{2}|\mathcal E|$, at least $R_s/2$ children remain. 

Let $F'$ be the resulting set after going through all $S$ blocks. By the uniform branching of $F$,
\begin{equation}\label{eq:markov-retained-cardinality}
    |F'|\ge2^{-S}|F|.
\end{equation}
Theorem \ref{thm:inverse-thm}(i), \eqref{eq:markov-retained-cardinality}, and Young's inequality yield
\begin{equation}\label{eq:markov-pruned-lower}
    \|\widetilde\mu^{(M)}|_{F'}\|_q^q
    \ge2^{-S-q-q\delta M}\|\theta^{(M)}*\mu_x^{(M)}\|_q^q.
\end{equation}
For $s\in\{0,\ldots,S-1\}\setminus\mathcal S_1$ and $Q\in\D_{sL}(F')$, let $\mathcal E=\D_{(s+1)L}(F'\cap Q)$. Every member of $\mathcal E$ meets $\frac12 Q$ by Theorem \ref{thm:inverse-thm}(ii).
For $\omega=x_0^\infty$, write
\[
    g_\omega=\tau_{t+\pi_{T_{sL+3}-1}(\omega)}S_{\lambda(x_0^{T_{sL+3}})},
    \qquad
    \Omega=\left\{\omega:(g_\omega\mu_{x_{T_{sL+3}}})\left(\cup \mathcal E\right)>0\right\}.
\]
For $\omega\in\Omega$, the component $g_\omega\mu_{x_{T_{sL+3}}}$ is supported in $Q$. Thus \eqref{eq:markov-stopped-decomposition} yields $P_x^\infty(\Omega)\le\widetilde{\mu}(Q)$. The cubes $g_\omega^{-1}J$ have side lengths $\sim 2^{-L}$. Lemmas \ref{lem:partition-comparison}(i) and \ref{lem:Tq} show that, for $L\gg 1$,
\[
\sum_{J\in\mathcal E}(g_\omega\mu_{x_{T_{sL+3}}})(J)^q
\lesssim_q2^{-(T_{\mathcal X}(q)-\delta/2)L}.
\]
Integrating this inequality over $\Omega$
and applying Jensen's inequality yields
\begin{equation}\label{eq:markov-local-centered-Lq}
    \sum_{J\in\mathcal E}\widetilde{\mu}(J)^q
    \le 2^{-(T_{\mathcal X}(q)-\delta)L}\widetilde{\mu}(Q)^q.
\end{equation}
Iterating \eqref{eq:markov-local-centered-Lq} for $s\notin\mathcal S_1$ and \eqref{eq:markov-local-unsaturation-Lq} for $s\in\mathcal S_1$ over the $S$ blocks, and using \eqref{eq:markov-many-saturated-blocks}, we obtain
\begin{equation}\label{eq:markov-pruned-upper}
    \|\widetilde\mu^{(M)}|_{F'}\|_q^q
    \le2^{-T_{\mathcal X}(q)M+\delta M-\eta|\mathcal S_1|L}
    \le2^{-(T_{\mathcal X}(q)+\eta\sigma/(2d(q-1))-\delta)M}.
\end{equation}
On the other hand, the failure of \eqref{eq:markov-convolution-flattening-limit} and \eqref{eq:markov-pruned-lower}, together with our parameter choices, yield
\[
    \|\widetilde\mu^{(M)}|_{F'}\|_q^q
    >2^{-(T_{\mathcal X}(q)+\varepsilon+q\delta+L^{-1}+q/M)M}
    \stackrel{\eqref{eq:markov-delta-choice},\,\eqref{eq:markov-L-choice},\,\eqref{eq:markov-epsilon-choice}}{>}2^{-(T_{\mathcal X}(q)+\eta\sigma/(2d(q-1))-\delta)M}
\]
for large enough $M$. This is the desired contradiction.

The case of general $M$ follows by interpolating between $L\lfloor M/L\rfloor$ and $L\lceil M/L\rceil$, and absorbing the resulting constant by making $\varepsilon$ slightly smaller.
\end{proof}

\begin{cor}\label{cor:unsaturated-subexponential-flattening}
If a Markov self-similar model is $q$-unsaturated on lines and has sub-exponential endpoint growth, then it satisfies the $L^q$ flattening condition.
\end{cor}

\begin{proof}
The $L^q$ flattening condition follows from Lemma \ref{lem:subexponential-orbit-growth-flattening} and Theorem \ref{thm:markov-convolution-flattening}.
\end{proof}

\section{Diophantine and transversal systems}

In this section, we specialize to a class of Markov self-similar models that satisfy additional Diophantine and transversality conditions. We will later verify that the models underpinning Theorems \ref{thm:intro-slicing} and \ref{thm:intro-ssm-slices} satisfy these conditions.

In this section, given a Markov self-similar model $\mathcal{X}=(X,P,\pi,\lambda)$, we fix a metric $d$ on $X$ inherited from the compact metric space including $X$.

\begin{dfn} \label{dfn:diophantine}
    Fix $K> 0$. We say the Markov self-similar model is \emph{$K$-Diophantine} if there are infinitely many $N$ such that for every Markov orbit $(y_0,y_1,\ldots,y_N)$, we have
    \[
        d(y_i,y_j)\gtrsim N^{-K}
        \qquad\text{for all }0\le i<j\le N,
    \]
    where the implicit constant is independent of the orbit and of $N$.
\end{dfn}

This Diophantine condition is very concrete, but the proofs allow for a weaker version. We denote by $|Y|_\delta$ the $\delta$-covering number of a set $Y\subset X$, that is, the smallest number of closed balls of radius $\delta$ needed to cover $Y$.

\begin{dfn} \label{dfn:weakly-diophantine}
    Fix $c>0$. We say the Markov self-similar model is \emph{weakly $c$-Diophantine} if there are sequences $r_j\to 0$, $N_j\to\infty$ such that for every Markov orbit $(y_0,y_1,\ldots,y_{N_j})$ and every set $E\subset\{0,1,\ldots,N_j\}$, we have
    \[
        \left| \{ y_i : i\in E\} \right|_{r_j} \gtrsim \underline d_S(E)\cdot r_j^{-c},
    \]
where the implicit constant is independent of $j$, the orbit, and $E$.
\end{dfn}
If a model is $K$-Diophantine then it is weakly $(1/K)$-Diophantine. Indeed, for each $N$ in Definition \ref{dfn:diophantine}, choose $r\sim N^{-K}$ so that the states along every orbit are $3r$-separated. Each ball of radius $r$ contains at most one selected state, so the covering number is at least $|E|=(\text{density of $E$})N\gtrsim \underline d_S(E) r^{-1/K}$.

On a first reading of the following next two definitions, the reader may ignore the partitions and simply take $C=\supp(P_x^{[m]})$. The flexibility provided by the partitions will, however, be important in our applications.
\begin{dfn} \label{dfn:transversal}
Fix $L,c\ge0$. A Markov self-similar model is \emph{$(L,c)$-transversal} if, for each $x\in X$ and $m\in\N$, there are finite partitions $\mathcal C_{x,m}$ of $\supp(P_x^{[m]})$ satisfying
\begin{equation}\label{eq:separation-class-growth}
\lim_{m\to\infty}\frac{1}{m}\sup_{x\in X}\log|\mathcal C_{x,m}|= 0,
\end{equation}
such that, for every $m\in\N$ and $0<r\le1$, the set
\[
E_{m,r}=\left\{x\in X:\begin{array}{l}
\text{there are }C\in\mathcal C_{x,m}\text{ and distinct }\mathbf{x},\mathbf{y}\in C\\
\text{with }|\widetilde \pi(\mathbf{x})-\widetilde \pi(\mathbf{y})|\le r
\end{array}\right\}
\]
can be covered by $\lesssim2^{Lm}r^{-c}$ closed balls of radius $r$, with an implicit constant independent of $m,r$.
\end{dfn}

\begin{dfn} \label{dfn:uniform-exponential-separation}
A Markov self-similar model is \emph{uniformly exponentially separated} if there are $R_0>1$, an infinite set $S\subset\N$ and finite partitions $\mathcal C_{x,m}$ of $\supp(P_x^{[m]})$ for every $x\in X$ and $m\in\N$, satisfying \eqref{eq:separation-class-growth}, such that, for every $m\in S$, $x\in X$, $C\in\mathcal C_{x,m}$ and distinct $\mathbf{x},\mathbf{y}\in C$,
\[
|\widetilde\pi(\mathbf{x})-\widetilde\pi(\mathbf{y})|\ge2^{-R_0m}.
\]
\end{dfn}

The following is our main technical result on $L^q$ dimensions of Markov self-similar measures under Diophantine and transversality conditions.
\begin{thm} \label{thm:diophantine-transversal-Lq}
    Consider a Markov self-similar model $\mathcal{X}$ that satisfies the $L^q$ flattening condition for some $q>1$ and such that all transition laws are purely atomic. Assume, additionally, that either $\mathcal{X}$ is uniformly exponentially separated, or it satisfies both the weakly $c$-Diophantine and the $(L,c')$-transversal condition for some $c>c'$. Then,
   \[
      T(q) \ge \liminf_{m\to\infty} \inf_{x\in X} \frac{-\log\|P_x^{[m]}\|_q^q}{m}.
   \]
\end{thm}
\begin{proof}
We begin with a counting observation used in both cases. Let $\mathcal C_{x,m}$ be the partitions from Definition \ref{dfn:transversal} or \ref{dfn:uniform-exponential-separation}. If $r\ge 2^{-k}$ and
\[
    |\widetilde\pi(\mathbf{x})-\widetilde\pi(\mathbf{y})|\ge r\quad\text{for }\mathbf{x}\neq\mathbf{y} \in C, \, C\in\mathcal C_{x,m},
\]
then $|C\cap\widetilde\pi^{-1}(Q)|\lesssim_d1$ for every $C\in\mathcal C_{x,m}$ and $Q\in\D_k$. Consequently, H\"older's inequality and \eqref{eq:separation-class-growth} yield
\begin{equation}\label{eq:stopped-path-multiplicity}
    \|\mu_{x,[m]}^{(k)}\|_q^q
    \lesssim_{d,q}|\mathcal C_{x,m}|^{q-1}\|P_x^{[m]}\|_q^q
    \le2^{o(m)}\|P_x^{[m]}\|_q^q,
\end{equation}
with $o(m)$ as $m\to\infty$ uniform in $x,r,k$.

Let $\e(\sigma)\ (\sigma>0)$ be the flattening rate given by the $L^q$ flattening condition for $\mathcal{X}$. For uniform exponential separation with parameters $R_0,S$, we may assume that $R_0\ge2$ is an integer. Fix $\sigma>0$ and $\e=\e(\sigma)\le\sigma$. For arbitrarily large $m\in S$, the definition of $T(q)$ ensures that some $x\in X$ satisfies $\phi_{(R_0+1)m}(x)\ge2^{-(R_0+1)m(T(q)+\e/2)}$. From Proposition \ref{prop:Lq-at-finer-scale} and \eqref{eq:stopped-path-multiplicity}, we then obtain
\[
2^{-m(T(q)+2R_0\sigma)}
\lesssim\|\mu_{x,[m]}^{((R_0+1)m)}\|_q^q
\le 2^{o(m)}\|P_x^{[m]}\|_q^q.
\]
Letting $m\to\infty$ and then $\sigma\to0$ settles the uniformly exponentially separated case.

Suppose now that the model is weakly $c$-Diophantine and $(L,c')$-transversal, with $c'<c$. Choose an integer $R\ge2$ with $L/(R+1)<(c-c')/2$, then fix $\sigma>0$ and $\e=\e(\sigma)\le\sigma$. For $r=r_j$ and $N=N_j$ from Definition \ref{dfn:weakly-diophantine}, set $m=\lceil\log(1/r)/(R+1)\rceil$. Then
\begin{equation}\label{eq:choice-of-R-stopped}
    2^{-(R+1)m}\le r,
    \qquad
    2^{Lm}r^{-c'}\lesssim r^{-c'-L/(R+1)}
    =o\bigl(r^{-(c+c')/2}\bigr).
\end{equation}
For large $j$, Proposition \ref{prop:good-step-n-orbit} provides a stopped scale-$(R+1)m$ orbit $(\tilde y_0,\ldots,\tilde y_N)$ for which
\[
\widetilde E=\bigl\{i\in\{0,\ldots,N\}:\phi_{(R+1)m}(\tilde y_i)>2^{-(R+1)m(T(q)+\e/2)}\bigr\}
\]
has Schnirelmann density at least $1/2$.

Choose a Markov orbit $(y_n)_{n\ge0}$ realizing $\tilde y_i=y_{t_i}$ at successive scale-$(R+1)m$ stopping times. Since $t_{i+1}-t_i\le (R+1)m/\alpha_{\min}+1$, the set $E=\{t_i:i\in\widetilde E,\ t_i\le N\}$ has $\underline d_S(E)\gtrsim1/(R+1)m$. By the weakly Diophantine condition,
\[
\bigl|\{y_i:i\in E\}\bigr|_r\gtrsim m^{-1}r^{-c}
\ge r^{-(c+c')/2}
\]
for large $j$, since $m\sim\log(1/r)$. Definition \ref{dfn:transversal} and \eqref{eq:choice-of-R-stopped} yield $|E_{m,r}|_r=o(r^{-(c+c')/2})$. Thus some $x\in\{y_i:i\in E\}\setminus E_{m,r}$ satisfies 
\[
\phi_{(R+1)m}(x)>2^{-(R+1)m(T(q)+\e/2)}, 
\]
and $|\widetilde\pi(\mathbf{x})-\widetilde\pi(\mathbf{y})|>r$ for distinct $\mathbf{x},\mathbf{y}$ in any $C\in\mathcal C_{x,m}$. Applying Proposition \ref{prop:Lq-at-finer-scale} and \eqref{eq:stopped-path-multiplicity}, we obtain
\[
2^{-m(T(q)+2R\sigma)}
\lesssim\|\mu_{x,[m]}^{((R+1)m)}\|_q^q
\le 2^{o(m)}\|P_x^{[m]}\|_q^q.
\]
Letting $j\to\infty$ and then $\sigma\to0$, we get the claim.
\end{proof}

\section{Intersections of self-similar sets}\label{sec:embedding_conjecture}

\subsection{Main result and discussion of the Diophantine condition}

\label{sec:diophantine-convolution}

Consider two linear iterated function systems $\Phi=\{\lambda_i x+a_i\}_{i\in \mathcal{I}}$ and $\Psi=\{\lambda'_jx+b_j\}_{j\in \mathcal{J}}$ in $\R$. 

We assume strong exponential separation: there is $R\ge1$ such that, for every $n\ge1$, compositions $f,g$ indexed by distinct words of length $n$ from the same IFS satisfy $|f(0)-g(0)|\ge2^{-Rn}$. This condition follows from the open set condition after translating each IFS so that $0$ lies in one of its feasible open sets. It also holds if all contraction ratios and translations are algebraic and distinct words of the same length have distinct images of $0$; see \cite[Proof of Theorem 1.5]{Hoc14}.

We write $\alpha_i = \log(1/\lambda_i)$ and $\beta_j= \log(1/\lambda'_j)$ for simplicity.

\begin{dfn} \label{dfn:diophantine-convolution}
    Fix $K>0$. We say that $(\Phi,\Psi)$ are \emph{$K$-Diophantine} if there are $c_0>0$ and infinitely many $N$ such that for all integers $0\le n_i, m_j\le N\ (i\in\mathcal{I},j\in\mathcal{J})$, with at least one of them non-zero, we have
    \[
        \left| \sum_{i\in\mathcal{I}} n_i\alpha_i -\sum_{j\in\mathcal{J}} m_j\beta_j \right| \ge c_0N^{-K}.
    \]
\end{dfn}
This can be seen as a quantitative form of \eqref{eq:no-exact-overlaps-slicing}. Note that the condition depends only on the contraction ratios of the maps, and not on the translations. The next theorem is a stronger form of Theorem \ref{thm:intro-slicing} in the introduction.
\begin{thm} \label{thm:convolution-ssm}
    Let $\Phi=\{\lambda_i x+a_i\}_{i\in \mathcal{I}}$ and $\Psi=\{\lambda'_jx+b_j\}_{j\in \mathcal{J}}$ be two linear iterated function systems in $\R$ satisfying the strong exponential separation condition. Assume that $(\Phi,\Psi)$ are $K$-Diophantine for some $K>0$. 
    
    Let $\nu,\nu'$ be self-similar measures for $\Phi, \Psi$ and probability vectors $p=(p_i)_{i\in\mathcal{I}}, p'=(p'_j)_{j\in\mathcal{J}}$, respectively. Then, for all $q>1$,
    \[
        D_q(\nu*\nu') = \min\left\{\frac{\tau_{\Phi,p}(q)+\tau_{\Psi,p'}(q)}{q-1},1\right\},
    \]
    where $\tau_{\Phi,p}(q)$ is the unique real number satisfying $\sum_{i\in\mathcal{I}}p_i\lambda_i^{-\tau_{\Phi,p}(q)}=1$ and $\tau_{\Psi,p'}(q)$ is defined similarly.
\end{thm}

According to Remark \ref{rem:slices}, we get:
\begin{cor}
    Under the assumptions of Theorem \ref{thm:convolution-ssm}, let $A$ and $B$ be the self-similar sets for $\Phi$ and $\Psi$, respectively. Then
    \[
        \hdim\bigl(A\cap g(B)\bigr) \le \max\left\{0,\hdim(A)+\hdim(B)-1\right\}
    \]
    for all affine maps $g:\R\to\R$.
\end{cor}

The Diophantine condition in Definition \ref{dfn:diophantine-convolution} is very similar to condition $(d)$ in \cite[Definition 6.2]{AHW24}. More precisely, it is shown in \cite{AHW24} that if $\Phi$ is homogeneous, both $\Phi$ and $\Psi$ satisfy the strong separation condition, and $(\Phi,\Psi)$ are Diophantine in the sense of Definition \ref{dfn:diophantine-convolution}, then the attractor of $\Phi$ does not embed affinely into the attractor of $\Psi$. Our corollary strengthens this in several ways. Most importantly, it provides not only a non-embedding result, but also a ``Furstenberg slicing'' bound on the intersection. Secondly, it allows for non-homogeneous $\Phi$. Lastly, we substantially relax the strong separation condition to the strong exponential separation condition.

A necessary condition for Definition \ref{dfn:diophantine-convolution} is
\begin{equation}\label{eq:diophantine-no-product-relation}
    \prod_{i\in\mathcal I}\lambda_i^{n_i}
    \ne\prod_{j\in\mathcal J}(\lambda_j')^{m_j}
    \qquad
    \bigl(\mathbf n\in\N_0^{\mathcal I},\ \mathbf m\in\N_0^{\mathcal J},\ (\mathbf n,\mathbf m)\ne0\bigr).
\end{equation}
When all contraction ratios are algebraic, \eqref{eq:diophantine-no-product-relation} is also sufficient by Baker's theorem \cite[Theorem 6.5 and \S6.2]{AHW24}: every non-zero linear form in Definition \ref{dfn:diophantine-convolution} is bounded below by a fixed negative power of $N$. In particular, it suffices that $\spn_{\Q}\{\alpha_i\}_i\cap\spn_{\Q}\{\beta_j\}_j=\{0\}$; linear relations within either family are allowed. For example, both families may consist of positive integer powers of $1/2$ and $1/3$, respectively.

Condition \eqref{eq:diophantine-no-product-relation} is also sufficient when all $\alpha_i,\beta_j$ are a common positive multiple of algebraic numbers, as noted in \cite[\S1]{AHW24}: the algebraic norm bounds the non-zero linear forms below by a fixed negative power of $N$. In both cases, the Diophantine inequality holds for all sufficiently large $N$.

There is also an explicit case requiring no algebraicity. Suppose $\lambda_i=r^{k_i}$ and $\lambda_j'=s^{\ell_j}$ for $r,s\in(0,1)$ and positive integers $k_i,\ell_j$. Then $(\Phi,\Psi)$ are Diophantine if and only if $\log r/\log s\notin\Q$. In fact, $K=1$ works in the irrational case:
for $\alpha=\log r/\log s\notin\Q$, if there are integers $0<n_N,m_N\le N$ such that $|n_N\alpha-m_N|<10^{-1}N^{-1}$ for every large $N\in\N$, then, by estimating $|m_N/n_N-m_{N+1}/n_{N+1}|$ for every $N$, $m_N/n_N$ turns out to be constant and equal to $\alpha$ and we get a contradiction.
If $\log r/\log s\in\Q$, then \eqref{eq:diophantine-no-product-relation} fails. Thus this case includes two non-homogeneous IFSs whose logarithmic base ratio is Liouville. 

When $\Phi$ is homogeneous, write $\lambda$ for its common contraction ratio. The following criteria can be quoted directly from \cite{AHW24}. If $\lambda$ and all $\lambda'_j$ are algebraic and $\log(1/\lambda)\notin\spn_{\Q}\{\beta_j\}_j$, then $(\Phi,\Psi)$ are Diophantine \cite[\S6.2]{AHW24}. For fixed $\Psi$, the exceptional set of $\lambda\in(0,1)$ has packing dimension zero \cite[Lemma 6.7]{AHW24}; for fixed $\lambda$, the exceptional set of $(\lambda'_j)_j\in(0,1)^{\mathcal J}$ has packing dimension at most $|\mathcal J|-1$ \cite[Lemma 6.8]{AHW24}. Here ``exceptional'' means that Definition \ref{dfn:diophantine-convolution} fails for every $K>0$, and $\Psi$ need not be homogeneous.

The covering argument in \cite[Lemmas 6.7--6.8]{AHW24} extends to two non-homogeneous IFSs as follows. Let $\mathcal E\subset(0,\infty)^{\mathcal I}\times(0,\infty)^{\mathcal J}$ be the set of logarithmic contraction vectors $(\alpha,\beta)=((\alpha_i)_i,(\beta_j)_j)$ for which Definition \ref{dfn:diophantine-convolution} fails for every $K>0$. Then
\begin{equation}\label{eq:diophantine-exceptional-dimensions}
    \begin{aligned}
        \pdim\mathcal E&=|\mathcal I|+|\mathcal J|-1,\\
        \pdim\{\alpha:(\alpha,\beta)\in\mathcal E\}&=|\mathcal I|-1
            &&\text{for every fixed }\beta,\\
        \pdim\{\beta:(\alpha,\beta)\in\mathcal E\}&=|\mathcal J|-1
            &&\text{for every fixed }\alpha.
    \end{aligned}
\end{equation}
For \eqref{eq:diophantine-exceptional-dimensions}, the argument of \cite[Lemmas 6.7--6.8]{AHW24} applies to the $N^{-K}$-neighborhoods of the $O(N^{|\mathcal I|+|\mathcal J|})$ hyperplanes $\sum_in_i\alpha_i=\sum_jm_j\beta_j$ in each compact box. With $\beta$ fixed, use affine hyperplanes in $\R^{\mathcal I}$; terms with $\mathbf n=0$ are excluded for large $N$ by positivity of $\beta$. Interchanging $\alpha,\beta$ yields the bound for $\pdim\{\beta:(\alpha,\beta)\in\mathcal E\}$. The relations $\alpha_i=\beta_j$ show sharpness.

There is also a one-parameter version, as in \cite[Lemma 6.7]{AHW24}. For fixed $\alpha,\beta$, replace the contractions of $\Phi$ by $(\lambda_i^s)_i$, $s>0$. Then
\[
    \pdim\{s>0:(s\alpha,\beta)\in\mathcal E\}=0,
\]
so almost every $s>0$ satisfies the Diophantine condition, even when both IFSs are non-homogeneous. The same covering argument applies because $|s\sum_in_i\alpha_i-\sum_jm_j\beta_j|<N^{-K}$ defines $O(N^{|\mathcal I|+|\mathcal J|})$ intervals of length $O(N^{-K})$ on compact intervals of positive $s$.

For comparison, under strong separation, \cite[Theorems 1.3--1.4]{AHW24} obtain stronger conclusions about non-embedding by selecting homogeneous subsystems. For fixed $\Psi$, the exceptional vectors of contractions of $\Phi$ for their embedding conjecture have packing dimension zero. For fixed algebraic contractions of $\Phi$, the exceptional vectors for $\Psi$ have packing dimension at most $|\mathcal J|-\dim_{\Q}\spn_{\Q}\{\alpha_i\}_i$, and the exceptional set is empty if this rank exceeds $|\mathcal J|$. These are bounds for exceptions to non-embedding, rather than to Definition \ref{dfn:diophantine-convolution}; the latter can fail whenever $\alpha_i=\beta_j$ for even one pair $i,j$.

\subsection{The convolution model and its basic properties}

Let $f_i(x)=\lambda_i x+a_i$ and $g_j(x)=\lambda'_j x+b_j$ for $i\in\mathcal I$ and $j\in\mathcal J$. Let
\[
\nu=\sum_{i\in\mathcal I} p_i f_i\nu,
\qquad
\nu'=\sum_{j\in\mathcal J} p'_j  g_j\nu'
\]
be the given self-similar measures. Assume that
\[
\alpha_i<\beta_j
\qquad\text{for all }i\in\mathcal I,\ j\in\mathcal J.
\]
This can be arranged by replacing $\Psi$ by a suitable iterate, without changing the self-similar measure $\nu'$.

Combining the constructions of \S\ref{S:Markov-model-self-similar}, \S\ref{S:Markov-model-products}, and \S\ref{sec:convolution-model}, we obtain the following explicit model. Set
\[
X=\{(u,i,j): i\in\mathcal I,\ j\in\mathcal J,\ 0\le u<\beta_j\}.
\]
Define a Markov kernel $P$ on $X$ by
\[
P_{(u,i,j)}
=
\begin{cases}
\displaystyle \sum_{i'\in\mathcal I} p_{i'}\delta_{(u+\alpha_i,i',j)},
& \text{if }u+\alpha_i<\beta_j,\\
\displaystyle \sum_{i'\in\mathcal I}\sum_{j'\in\mathcal J} p_{i'}p'_{j'}\delta_{(u+\alpha_i-\beta_j,i',j')},
& \text{if }u+\alpha_i\ge \beta_j.
\end{cases}
\]
We also define
\[
\pi(u,i,j)
=
\begin{cases}
a_i,
& \text{if }u+\alpha_i<\beta_j,\\
a_i+2^u b_j,
& \text{if }u+\alpha_i\ge \beta_j,
\end{cases}
\qquad
\lambda(u,i,j)=\lambda_i.
\]
Then $\mathcal X=(X,P,\pi,\lambda)$ is a Markov self-similar model. The measures generated by $\mathcal X$ satisfy
\[
\mu_{(u,i,j)}=(f_i\nu)*S_{2^u}(g_j\nu').
\]
In particular,
\[
\mu_{(0,i,j)}=(f_i\nu)*(g_j\nu')
\qquad\text{for all }i\in\mathcal I,\ j\in\mathcal J.
\]
Consequently,
\[
\nu*\nu'=\sum_{i\in\mathcal I}\sum_{j\in\mathcal J} p_i p'_j\,\mu_{(0,i,j)}.
\]

For $\mathbf{i}=i_1\cdots i_n$, write $f_{\mathbf{i}}=f_{i_1}\circ\cdots\circ f_{i_n}$, $\lambda_{\mathbf{i}}=\prod_{r=1}^n\lambda_{i_r}$, $p_{\mathbf{i}}=\prod_{r=1}^n p_{i_r}$ and $\alpha_{\mathbf{i}}=\sum_{r=1}^n\alpha_{i_r}$, and define $g_{\mathbf{j}}$, $\lambda_{\mathbf{j}}'$, $p_{\mathbf{j}}'$ and $\beta_{\mathbf{j}}$ analogously. Write $\mathbf{i}^-$ for the word obtained by deleting the last symbol, and set $\lambda_\varnothing=\lambda_\varnothing'=1$.

\begin{lem}\label{lem:conv-model-transversal}
    If $(\Phi,\Psi)$ satisfy strong exponential separation, then there is $L=L(\Phi,\Psi)>0$ such that $\mathcal X$ is $(L,0)$-transversal in the sense of Definition \ref{dfn:transversal}.
\end{lem}
\begin{proof}
    For a stopped path $x_0^n$ starting at $x=(u,i,j)$, write $x_k=(u_k,i_k,j_k)\ (0\leq k\leq N)$ and $\mathbf{i}=(i_0,\ldots,i_n)$ and let $\mathbf{j}$ be the word formed by the symbols $j_k$, in increasing order of $k$, for which $0\le k\le n$ and $u_k+\alpha_{i_k}\ge\beta_{j_k}$. Then $\pi_n(x_0^n)=f_{\mathbf{i}}(0)+2^u g_{\mathbf{j}}(0)$, with $g_\varnothing=\operatorname{Id}$. The words $\mathbf{i},\mathbf{j}$ and $j_n$ determine the path by the transition rule for $P$ starting from $x$. Let $\mathcal C_{x,m}$ be the partition by $(|\mathbf{i}|,|\mathbf{j}|,j_n)$. Both word lengths are $O(m+1)$, so $|\mathcal C_{x,m}|=O((m+1)^2)$ uniformly in $x$, verifying \eqref{eq:separation-class-growth}.

    Fix $0<r\le1$. For two distinct paths in one class, write the difference of their images under $\widetilde\pi$ as $a+2^u b$. Strong exponential separation ensures that $(a,b)\ne(0,0)$. Increasing $R$ to account for the word lengths, each nonzero $a,b$ has magnitude at least $2^{-Rm}$. The set 
    \[
    \{u\in[0,\beta_j):|a+2^u b|\le r\}
    \] 
    is covered by $O(2^{Rm})$ balls of radius $r$: for $b\ne0$, use monotonicity and the derivative bound $|(2^u)'b|\gtrsim2^{-Rm}$; for $b=0$, the set is empty unless $r\ge2^{-Rm}$. There are $2^{O(m)}$ possible word pairs across all classes, so $|E_{m,r}|_r\lesssim2^{Lm}$ for some $L=L(\Phi,\Psi)$.
\end{proof}

We write $\tau_{\Phi}(q)=\tau_{\Phi,p}(q)$ and $\tau_{\Psi}(q)=\tau_{\Psi,p'}(q)$ for simplicity. They satisfy $\sum_{i\in\mathcal I} p_i^q \lambda_i^{-\tau_{\Phi}(q)}=\sum_{j\in\mathcal{J}}(p'_j)^q(\lambda'_j)^{-\tau_{\Psi}(q)}=1$.
For every $C>0$, an infinite word has $O_C(1)$ prefixes with $|\log\lambda_{\mathbf{i}}^{-1}-t|\le C$ and exactly one prefix with $\lambda_{\mathbf{i}^-}\ge2^{-t}>\lambda_{\mathbf{i}}$. Summing the Bernoulli cylinder weights $p_{\mathbf{i}}^q\lambda_{\mathbf{i}}^{-\tau_\Phi(q)}$ yields, uniformly in $t\ge0$,
\begin{equation}  \label{eq:symbolic-window-sum}
    \sum_{|\log\lambda_{\mathbf{i}}^{-1}-t|\le C}p_{\mathbf{i}}^q
    \lesssim_C 2^{-t\tau_\Phi(q)}, \qquad  \sum_{\lambda_{\mathbf{i}^-}\ge2^{-t}>\lambda_{\mathbf{i}}}p_{\mathbf{i}}^q
    \sim_q 2^{-t\tau_\Phi(q)}.
\end{equation}
The same relations hold for $\Psi$.

\begin{lem} \label{lem:conv-model-symbolic-Lq}
    For the convolution model $\mathcal X$, for every $m\in\N$ and every $x\in X$, 
    \[
       \|P_x^{[m]}\|_q^q \lesssim 2^{-m(\tau_{\Phi}(q)+\tau_{\Psi}(q))}.
    \]
    The implicit constant depends only on $\Phi,\Psi, q$ and not on $x$ or $m$.
\end{lem}
\begin{proof}
    Fix $x=(u,i,j)$ and take $m$ sufficiently large. Let $\mathbf{i}=i_1\cdots i_{T_m-1}$, and let $\mathbf{j}$ be the word formed by the symbols $j_k$, in increasing order of $k$, for which $0\le k<T_m$ and $u_k+\alpha_{i_k}\ge\beta_{j_k}$, with the first letter omitted. By the stopping rule, $|\alpha_{\mathbf{i}}-m|,|\beta_{\mathbf{j}}-m|\le C$ for some $C$ depending only on $\Phi,\Psi$. These words and $i_{T_m},j_{T_m}$ determine the stopped path. Thus
    \[
        \|P_x^{[m]}\|_q^q
        \le \sum_{|\alpha_{\mathbf{i}}-m|\le C}p_{\mathbf{i}}^q
        \sum_{|\beta_{\mathbf{j}}-m|\le C}(p_{\mathbf{j}}')^q
        \stackrel{\eqref{eq:symbolic-window-sum}}{\lesssim} 2^{-m(\tau_\Phi(q)+\tau_\Psi(q))}.
        \qedhere
    \]
\end{proof}

\subsection{Proof of Theorem \ref{thm:convolution-ssm}}

We now complete the proof of Theorem \ref{thm:convolution-ssm}.  Covering by cylinders of size $\sim 2^{-m}$ as in e.g. \cite[(6.4)]{Shm19}, we obtain $\|\mu_x^{(m)}\|_q^q\gtrsim_q2^{-m(\tau_\Phi(q)+\tau_\Psi(q))}$. Also, $\|\mu_x^{(m)}\|_q^q\gtrsim_q2^{-m(q-1)}$ by Hölder's inequality. Hence
\[
    T_{\mathcal X}(q)\le\min\{\tau_\Phi(q)+\tau_\Psi(q),q-1\}.
\]
Suppose first that $T_{\mathcal X}(q)<q-1$.
Let $(u_{n},i_{n},j_{n})_{n=0}^N$ be a Markov orbit in $X$. Then, for each $0\le k<\ell\le N$, we have
\[
u_\ell-u_k = \sum_{i\in\mathcal I} n_i \alpha_i - \sum_{j\in\mathcal J} m_j \beta_j,
\]
for some non-negative integers $n_i, m_j$, all bounded by $N$. Definition \ref{dfn:diophantine-convolution} applies, and the weakly $(1/K)$-Diophantine condition follows from the observation after Definition \ref{dfn:weakly-diophantine}.
The model satisfies the $(L,0)$-transversal hypothesis of Theorem \ref{thm:diophantine-transversal-Lq} thanks to Lemma \ref{lem:conv-model-transversal}.

For an orbit started at $x=(u,i_0,j_0)$, write $x_n=(u_n,i_n,j_n)$. The symbol counts
\[
\bigl(\#\{0\le n<T_m:i_n=i\}\bigr)_{i\in\mathcal I},
\qquad
\bigl(\#\{0\le n<T_m:j_n=j,\ u_n+\alpha_{i_n}\ge\beta_{j_n}\}\bigr)_{j\in\mathcal J}
\]
and $i_{T_m},j_{T_m}$ determine both $\lambda(x_0^{T_m})$ and $x_{T_m}$. Since $T_m\lesssim m+1$, there are $O((m+1)^{|\mathcal I|+|\mathcal J|})$ possible pairs in $\mathcal O_m(x)$. Thus the convolution model has sub-exponential endpoint growth as in Definition \ref{dfn:subexponential-orbit-growth}.
The $L^q$ flattening condition follows from Corollary \ref{cor:unsaturated-subexponential-flattening}. Together with Theorem \ref{thm:diophantine-transversal-Lq} and Lemma \ref{lem:conv-model-symbolic-Lq}, we conclude that $T_{\mathcal X}(q)\ge\tau_\Phi(q)+\tau_\Psi(q)$.

In either case, $T_{\mathcal X}(q)=\min\{\tau_\Phi(q)+\tau_\Psi(q),q-1\}$. The identity $\nu*\nu'=\sum_{i,j}p_i p'_j\mu_{(0,i,j)}$, Jensen's inequality, and Lemma \ref{lem:Tq} show that $\underline\tau_{\nu*\nu'}(q)\ge T_{\mathcal X}(q)=\min\{\tau_\Phi(q)+\tau_\Psi(q),q-1\}$. Since the upper bound is trivial,
we obtain $\tau_{\nu*\nu'}(q)=\min\{\tau_\Phi(q)+\tau_\Psi(q),q-1\}$.
This establishes the claimed formula for $D_q(\nu*\nu')$.

\subsection{A remark on the higher dimensional case}

By combining the present approach with the strategy of \cite[Theorem 1.15]{CS24} on slices of products of homogeneous self-similar measures, one can obtain higher dimensional versions of Theorem \ref{thm:convolution-ssm}, concerning projections of products of self-similar measures. We focused on the planar case for simplicity, and leave it to the reader to explore the details of the higher dimensional generalization.

\section{Planar self-similar sets and their projections: proof of Theorem \ref{thm:intro-ssm-slices}}\label{sec:planar-self-similar-projections}

\subsection{Introduction and main result}

We now apply the Markov formalism to non-homogeneous planar self-similar measures. Let $\mathcal I$ be a finite alphabet with at least two elements, and consider an orientation-preserving planar IFS
\[
    \Phi=\{f_i(z)=\lambda_iR_{\theta_i}z+t_i:i\in\mathcal I\},
\]
where $R_\theta$ denotes rotation through the angle $2\pi\theta$, $\theta_i\in\R/\Z$, $t_i\in\R^2$, and $0<\lambda_i<1$. Let $p=(p_i)_{i\in\mathcal I}$ be a positive probability vector and let $\mu$ be the corresponding self-similar measure. We use the model $\mathcal X$ from \S\ref{S:Markov-model-self-similar}, with
\[
    X=\SO(2)\times\mathcal I,
    \qquad
    P_{(O,i)}=\sum_{j\in\mathcal I}p_j\delta_{(OR_{\theta_i},j)},
    \qquad
    \pi(O,i)=Ot_i,
    \qquad
    \lambda(O,i)=\lambda_i,
\]
Recall that $\mu_{(O,i)}=O(f_i\mu)$.

We impose a Diophantine condition on the rotations: there are $c_0,K>0$ and infinitely many $N\in\N$ such that
\begin{equation}\label{eq:planar-angular-diophantine}
    \left\|\sum_{i\in\mathcal I}n_i\theta_i\right\|_{\R/\Z}
    \ge c_0N^{-K}
\end{equation}
whenever $n_i\in\N_0$, $0\le n_i\le N$ and at least one of them is non-zero. In particular, the group generated by the rotations $R_{\theta_i}$ is dense in $\SO(2)$.

For a word $\mathbf{i}=i_1\cdots i_n$, write $f_{\mathbf{i}}=f_{i_1}\circ\cdots\circ f_{i_n}$, $\lambda_{\mathbf{i}}=\lambda_{i_1}\cdots\lambda_{i_n}$, $p_{\mathbf{i}}=p_{i_1}\cdots p_{i_n}$, and let $\mathbf{i}^-=i_1\cdots i_{n-1}$.

We assume the strong exponential separation condition from \S\ref{sec:embedding_conjecture}: there is $R\ge1$ such that, for every $n\ge1$ and distinct $\mathbf{i},\mathbf{j}\in\mathcal I^n$,
\begin{equation}\label{eq:planar-strong-exponential-separation}
    |f_{\mathbf{i}}(0)-f_{\mathbf{j}}(0)|\ge2^{-Rn}.
\end{equation}
For $q>1$, let $\tau_\Phi(q)$ be the unique positive number satisfying
\begin{equation}\label{eq:planar-symbolic-spectrum}
    \sum_{i\in\mathcal I}p_i^q\lambda_i^{-\tau_\Phi(q)}=1.
\end{equation}
As before, this is the symbolic $L^q$ spectrum of $\mu$, and equals the $L^q$ spectrum under the open set condition and is always an upper bound for the $L^q$ spectrum.

Since strong exponential separation follows from the open set condition after a translation, the next theorem contains Theorem \ref{thm:intro-ssm-slices} from the introduction as a special case.
\begin{thm}\label{thm:planar-self-similar-Lq}
Assume \eqref{eq:planar-angular-diophantine} and the strong exponential separation condition \eqref{eq:planar-strong-exponential-separation}. Then, for every $q>1$ and every one-dimensional linear subspace $V\subset\R^2$,
\begin{equation}\label{eq:planar-model-spectrum}
    T_{\mathcal X}(q)=\min\{\tau_\Phi(q),2(q-1)\}
\end{equation}
and
\begin{equation*}
T_{\Pi_V\mathcal{X}}(q)=\min\left\{\tau_{\Phi}(q),q-1\right\}.
\end{equation*}
Consequently,
\[
    D_q(\mu)=\min\left\{\frac{\tau_\Phi(q)}{q-1},2\right\},
    \qquad
    D_q(\Pi_V\mu)=\min\left\{\frac{\tau_\Phi(q)}{q-1},1\right\}
\]
for every $V\in\mathbb{G}(2,1)$. The estimates are uniform in $V$: for every $\rho>0$ and all sufficiently large $m$,
\[
    \sup_{V\in\mathbb{G}(2,1)}\|(\Pi_V\mu)^{(m)}\|_q^q
    \le2^{-m(\min\{\tau_\Phi(q),q-1\}-\rho)}.
\]
Consequently, writing $A=\supp(\mu)$,
\[
 \ubdim(A\cap\ell) \le \max\{ \hdim(A)-1,0\}\quad\text{for every affine line }\ell\subset\R^2. 
\]
\end{thm}

We record several cases of \eqref{eq:planar-angular-diophantine} which parallel those in \S\ref{sec:diophantine-convolution}. A necessary condition is
\begin{equation}\label{eq:planar-no-rotation-relation}
    \sum_{i\in\mathcal I}n_i\theta_i\ne0\quad\text{in }\R/\Z
    \qquad\text{for all }\mathbf n\in\N_0^{\mathcal I}\setminus\{0\}.
\end{equation}
As in \S\ref{sec:diophantine-convolution}, condition \eqref{eq:planar-no-rotation-relation} is also sufficient if the angles have algebraic representatives, or if all rotation matrices have algebraic entries; the latter case uses the Baker--W\"ustholz bound \cite{BW93}. In both cases, \eqref{eq:planar-angular-diophantine} holds for all sufficiently large $N$.

If $\theta_i=a_i\theta$ modulo $1$, with positive integers $a_i$ and $\theta\notin\Q$, then \eqref{eq:planar-angular-diophantine} holds with $K=1$ along infinitely many $N$, by the same continued-fraction argument as in \S\ref{sec:diophantine-convolution}. This includes a common irrational rotation, even when the angle is Liouville.

Let $\mathcal E\subset(\R/\Z)^{\mathcal I}$ consist of the angle vectors for which \eqref{eq:planar-angular-diophantine} fails for every $c_0,K>0$. The argument for \eqref{eq:diophantine-exceptional-dimensions}, following \cite[Lemmas 6.7--6.8]{AHW24}, yields
\begin{equation}\label{eq:planar-angular-exceptional-dimension}
    \pdim\mathcal E\le |\mathcal{I}|-1.
\end{equation}
The one-parameter argument in \S\ref{sec:embedding_conjecture} likewise shows that, for any fixed angle vector $(\theta_i)_i$,
\[
    \pdim\{s\in\R/\Z:(\theta_i+s)_{i\in\mathcal I}\in\mathcal E\}=0.
\]

The proof of Theorem \ref{thm:planar-self-similar-Lq} follows the project-and-induct strategy from \cite{CS24}: we first show that every one-dimensional projection of the planar self-similar measure has the expected $L^q$ dimension, see Proposition \ref{prop:planar-projection-spectrum}. In turn, this provides the required $q$-unsaturation property to handle the planar case. Similar to \cite{CS24}, the method extends to higher dimension under assumptions similar to \cite{CS24}, except that homogeneity can be replaced by commutativity of the orthogonal parts and a corresponding Diophantine condition. We do not pursue this here, since we hope the planar case is already of interest and the higher-dimensional case is more technical.

\subsection{Polynomial endpoint growth}

\label{S:planar-polynomial-orbit-growth}

Let
\[
    \mathcal W_i(m)=\{\mathbf{i}:i_1=i,\ \lambda_{\mathbf{i}^-}\ge2^{-m}>\lambda_{\mathbf{i}}\}.
\]
It follows from \eqref{eq:symbolic-window-sum} that
\begin{equation}\label{eq:planar-stopped-path-norm}
    \|P_{(O,i)}^{[m]}\|_q^q
    =\frac{1}{p_i^q}
      \sum_{\mathbf{i}\in\mathcal W_i(m)}p_{\mathbf{i}}^q
    \sim_q2^{-m\tau_\Phi(q)}.
\end{equation}

For an orbit started at $x=(O,i)$, write $\mathbf{i}=i_0\cdots i_{T_m-1}\in\mathcal W_i(m)$ and $j=i_{T_m}$. Its stopped atom is $Of_{\mathbf{i}}(t_j)=Of_{\mathbf{i}j}(0)$. Since $|\mathbf{i}|\lesssim m+1$, we may increase $R$ to obtain
\begin{equation}\label{eq:planar-equal-length-separation}
    |f_{\mathbf{i}}(t_j)-f_{\mathbf{j}}(t_\ell)|\ge2^{-Rm}
\end{equation}
whenever $(\mathbf{i},j)\ne(\mathbf{j},\ell)$ and $|\mathbf{i}|=|\mathbf{j}|$, with $\mathbf{i},\mathbf{j}\in\mathcal W_i(m)$. Take
\[
    \mathcal C_{(O,i),m}=\bigl\{\supp(P_{(O,i)}^{[m]})\cap X^{n+1}:n\ge1\bigr\}\setminus\{\varnothing\}.
\]
This partitions the stopped paths by their number of transitions, with $|\mathcal C_{(O,i),m}|=O(m+1)$. Since $O$ preserves distances, \eqref{eq:planar-equal-length-separation} verifies Definition \ref{dfn:uniform-exponential-separation} with these partitions and $S=\N$.

For $x=(O,i_0)$, the symbol counts $\bigl(\#\{0\le n<T_m:i_n=j\}\bigr)_{j\in\mathcal I}$ and $i_{T_m}$ determine both $\lambda(x_0^{T_m})$ and $x_{T_m}$, since planar rotations commute. Since $T_m\lesssim m+1$, there are $O((m+1)^{|\mathcal I|})$ possible pairs in $\mathcal O_m(x)$. Thus the planar model and each projected model have polynomial, and therefore sub-exponential, endpoint growth as in Definition \ref{dfn:subexponential-orbit-growth}.

\subsection{Linear projections}

Fix $V\in\mathbb{G}(2,1)$. Recall that $\Pi_V$ denotes orthogonal projection onto $V$.  The models $\Pi_V\mathcal X$ are isomorphic as $V$ varies: left multiplication in $\SO(2)$ changes $V$ to any prescribed line and commutes with every $R_{\theta_i}$. In particular, $T_{\Pi_V\mathcal X}(q)$ is independent of $V$.

\begin{prop}\label{prop:planar-projection-spectrum}
For every $V\in\mathbb{G}(2,1)$ and $q>1$,
\begin{equation}\label{eq:planar-projection-spectrum}
    T_{\Pi_V\mathcal X}(q)=\min\{\tau_\Phi(q),q-1\}.
\end{equation}
Also,
\begin{equation}\label{eq:planar-projection-measure-dimension}
D_q(\Pi_V\mu)=\min\left\{\frac{\tau_\Phi(q)}{q-1},1\right\}.
\end{equation}
Moreover, for every $\rho>0$ and all sufficiently large $m$,
\begin{equation}\label{eq:uniform-planar-projection-Lq}
    \sup_{V\in\mathbb{G}(2,1)}\|(\Pi_V\mu)^{(m)}\|_q^q
    \le2^{-m(\min\{\tau_\Phi(q),q-1\}-\rho)}.
\end{equation}
\end{prop}
\begin{proof}

From \eqref{eq:planar-stopped-path-norm} and H\"older's inequality we deduce the standard inequality
\begin{equation}\label{eq:planar-projection-upper}
    T_{\Pi_V\mathcal X}(q)\le\min\{\tau_\Phi(q),q-1\}.
\end{equation}
To establish \eqref{eq:planar-projection-spectrum}, we may assume that $T_{\Pi_V\mathcal X}(q)<q-1$. Thus, the model is $q$-unsaturated on lines.

For two states at times $s<t\le N$, the quotient of their rotation coordinates is $\prod_{n=s}^{t-1}R_{\theta_{i_n}}$. Thus \eqref{eq:planar-angular-diophantine} implies that $\Pi_V\mathcal X$ is $K$-Diophantine, and hence weakly $(1/K)$-Diophantine by the observation following Definition \ref{dfn:weakly-diophantine}.

For transversality, partition $\supp\bigl(P_{(O,i)}^{[m]}\bigr)$ by path length. Two distinct paths of equal length, represented by $(\mathbf{i},j)$ and $(\mathbf{j},\ell)$, are mapped by $\Pi_V\pi_{[m]}$ to $\Pi_V(Of_{\mathbf{i}}(t_j))$ and $\Pi_V(Of_{\mathbf{j}}(t_\ell))$, respectively. Let $a=f_{\mathbf{i}}(t_j)-f_{\mathbf{j}}(t_\ell)$; then $|a|\ge2^{-Rm}$ by \eqref{eq:planar-equal-length-separation}. For $0<r\le1$, the set $\{ O\in\text{SO}(2):|\Pi_V(Oa)|\le r\}$ can be covered by $O(2^{Rm})$ balls of radius $r$. There are $2^{O(m)}$ pairs of stopped paths, so the set in Definition \ref{dfn:transversal} satisfies
\begin{equation}\label{eq:planar-labelled-transversality}
    |E_{m,r}|_r\lesssim2^{Lm}
\end{equation}
for some $L$ independent of $V,m,r$. The $O(m+1)$ length classes satisfy \eqref{eq:separation-class-growth}, so $\Pi_V\mathcal X$ is $(L,0)$-transversal.

We saw in \S\ref{S:planar-polynomial-orbit-growth} that the model has sub-exponential endpoint growth. The $L^q$ flattening condition follows from Corollary \ref{cor:unsaturated-subexponential-flattening}. 

We have verified all assumptions of Theorem \ref{thm:diophantine-transversal-Lq}. By \eqref{eq:planar-stopped-path-norm}, we obtain $T_{\Pi_V\mathcal X}(q)\ge\tau_\Phi(q)$. Together with \eqref{eq:planar-projection-upper}, this gives \eqref{eq:planar-projection-spectrum}.

The identity $\Pi_V\mu=\sum_i p_i\mu^{\Pi_V\mathcal X}_{(\mathrm{Id},i)}$, Jensen's inequality, and Lemma \ref{lem:Tq} then yield the claim \eqref{eq:uniform-planar-projection-Lq}; uniformity in $V$ holds because the projected models $\Pi_V\mathcal X$ are isomorphic.

The sets $\Pi_V(f_{\mathbf{i}}(\supp\mu))$ with $\lambda_{\mathbf{i}^-}\ge2^{-m}>\lambda_{\mathbf{i}}$ have diameter $O(2^{-m})$, uniformly in $V$. It follows from \eqref{eq:symbolic-window-sum} that
\[
\|(\Pi_V\mu)^{(m)}\|_q^q
\gtrsim_{\Phi,p,q}\sum_{\lambda_{\mathbf{i}^-}\ge2^{-m}>\lambda_{\mathbf{i}}}p_{\mathbf{i}}^q
\sim_{\Phi,p,q}2^{-m\tau_\Phi(q)}.
\]
Together with the Hölder bound $\|(\Pi_V\mu)^{(m)}\|_q^q\gtrsim_{\Phi,p,q}2^{-m(q-1)}$ and \eqref{eq:uniform-planar-projection-Lq}, this yields \eqref{eq:planar-projection-measure-dimension}.
\end{proof}

\subsection{Proof of Theorem \ref{thm:planar-self-similar-Lq}}

We can now easily conclude the proof of Theorem \ref{thm:planar-self-similar-Lq}, and with it of Theorem \ref{thm:intro-ssm-slices}. It follows from \eqref{eq:planar-stopped-path-norm}  that
\[
    T_{\mathcal X}(q)\le\min\{\tau_\Phi(q),2(q-1)\}.
\]
Suppose $T_{\mathcal X}(q)<2(q-1)$. This and (\ref{eq:planar-projection-spectrum}) imply that the model $\mathcal{X}$ is $q$-unsaturated on lines.

We saw in \S\ref{S:planar-polynomial-orbit-growth} that the model has sub-exponential endpoint growth. The $L^q$ flattening condition follows from Corollary \ref{cor:unsaturated-subexponential-flattening}. The uniform exponential separation verified in \S\ref{S:planar-polynomial-orbit-growth} allows us to apply Theorem \ref{thm:diophantine-transversal-Lq}. Equation \eqref{eq:planar-stopped-path-norm} therefore yields $T_{\mathcal X}(q)\ge\tau_\Phi(q)$. This establishes \eqref{eq:planar-model-spectrum}.

Since $\mu_{(O,i)}=O(f_i\mu)$, it follows from Lemma \ref{lem:partition-comparison}(i) that $\|\mu_{(O,i)}^{(m)}\|_q^q\sim_{\Phi,p,q}\|\mu^{(m)}\|_q^q$, uniformly in $(O,i)$. Thus \eqref{eq:planar-model-spectrum} is the $L^q$ spectrum of $\mu$.

The projection formula and the uniform bound are \eqref{eq:planar-projection-measure-dimension} and \eqref{eq:uniform-planar-projection-Lq}. 
For the slicing bound, let $s$ satisfy $\sum_i\lambda_i^s=1$ and take $p_i=\lambda_i^s$  in the definitions of $\mu$ and $\tau_\Phi(q)$. Equation \eqref{eq:planar-model-spectrum} shows that $\hdim(A)=\min\{s,2\}$. If $s\ge2$, the claim follows from $\ubdim(A\cap\ell)\le1$. If $s<2$, then $\mu(B(x,r))\gtrsim r^s$ for $x\in A$ and $0<r<1$, so \eqref{eq:uniform-planar-projection-Lq} and \cite[Lemmas~1.7--1.8]{Shm19}, with $q\to\infty$, yield $\ubdim(A\cap\ell)\le\max\{s-1,0\}$.

\bibliographystyle{plain}
\bibliography{references}

\end{document}